\documentclass[a4paper,12pt,reqno]{amsart}
\usepackage[utf8]{inputenc}
\usepackage{amsmath}
\usepackage{amsthm,amssymb,amsfonts}
\usepackage{graphicx}
\usepackage{hyperref}
\usepackage{url}
\usepackage{mathrsfs}
\numberwithin{equation}{section}

\usepackage{geometry}
\usepackage{array,bm}
\usepackage{color}
\usepackage{verbatim} 

\newtheorem{theorem}{Theorem}[section]
\newtheorem{definition}[theorem]{Definition}
\newtheorem{corollary}[theorem]{Corollary}
\newtheorem{proposition}[theorem]{Proposition}
\newtheorem{remark}[theorem]{Remark}

\newtheorem{lemma}[theorem]{Lemma}

\newcommand{\Kn}{\mathcal{K}^n}
\newcommand{\Kno}{\mathcal{K}_o^n}
\newcommand{\Kne}{\mathcal{K}_e^n}
\newcommand{\Rn}{\mathbb{R}^n}
\newcommand{\Sn}{\mathbb{S}^{n-1}}
\newcommand{\bR}{\mathbb{R}}
\newcommand{\vol}{\operatorname{vol}}
\newcommand{\cHnmone}{\mathcal{H}^{n-1}}
\newcommand{\cH}{\mathcal{H}}

\newcommand{\bP}{\mathbb{P}}

\makeatletter
\@namedef{subjclassname@2020}{\textup{2020} Mathematics Subject Classification}
\makeatother

\begin{document}
	
	\title[MP for $(s,q)$-chord measures]{Minkowski problems for $(s,q)$-chord measures}

	\author{Jian-Yang Zhang}
	\address{Department of Mathematics, Shanghai University, Shanghai 200444, China}
	\email{3203927393@qq.com}
		\date{\today}
		\thanks{The authors' research was supported by NSFC grant no.~12322103.}
	
	\begin{abstract}
		Motivated by the work of Lutwak--Xi--Yang--Zhang on the chord measures,
		we introduce the $(s,q)$-chord measures arising from variations of a new family of geometric invariants, called
		the $(s,q)$-chord integrals. When $s=1$, the $(s,q)$-chord measures recover the chord measures. We study the corresponding Minkowski problems for the \((s,q)\)-chord measure.
	\end{abstract}
	
	\maketitle

	\section{Introduction}\label{sec:introduction}

Throughout the paper, $n\ge2$. The Euclidean inner product and norm on $\Rn$ are denoted by $x\cdot y$ and $|x|=(x\cdot x)^{1/2}$, respectively. We write $o$ for the origin, $\Sn$ for the unit sphere in $\Rn$, $B^n$ for the Euclidean unit ball centered at the origin, and $V$ for $n$-dimensional volume; set $\omega_n=V(B^n)$. For $u\in\Sn$, let
$u^\perp=\{x\in\Rn:x\cdot u=0\}$. A convex body $K$ is a compact convex subset of $\Rn$ with nonempty interior. The classes of convex bodies, convex bodies containing $o$ in their interiors, and origin-symmetric convex bodies (namely $K=-K$) are denoted by $\Kn$, $\Kno$, and $\Kne$, respectively. We write $\operatorname{SL}(n)$ for the group of determinant-one linear transformations of $\Rn$. Let $\cH^k$ be the $k$-dimensional Hausdorff measure for $k\ge0$. Unless otherwise stated, integration on $\Sn$ is with respect to $\cH^{n-1}$. For a finite Borel measure $\mu$ on $\Sn$, we write $|\mu|=\mu(\Sn)$.

The classical Minkowski problem asks for necessary and sufficient conditions on a finite Borel measure $\mu$ on $\Sn$ such that
\[
\mu=S(K,\cdot),
\]
where $K\subset\Rn$ is a convex body and $S(K,\cdot)$ is the surface area measure of $K$. The study of this problem was initiated by Minkowski
\cite{Minkowski1897,Minkowski1903}, who treated the cases
in which the prescribed measure is discrete or has a
continuous density. Later, Aleksandrov \cite{Aleksandrov1939} and, independently,
Fenchel--Jessen \cite{FenchelJessen1938} treated the case of
arbitrary finite Borel measures. From the analytic viewpoint, the problem amounts to solving a Monge--Amp\`ere equation. The regularity theory for its solutions is highly nontrivial and has motivated substantial work; see, for example, Nirenberg \cite{Nirenberg1953}, Cheng--Yau \cite{ChengYau1976}, Caffarelli \cite{Caffarelli1989,Caffarelli1990Localization,Caffarelli1990W2p}, and the references therein. We refer to Schneider's book \cite{schneiderbook} for the classical Brunn--Minkowski theory and its historical background.

Let $h_0,f$ be continuous functions on $\Sn$, and suppose that the Wulff shape $[h_0]$ is a convex body; see \eqref{defsec2-wulff-shape} for the definition of Wulff shape. The approach of Aleksandrov is variational and relies crucially on the formula
\begin{equation}\label{eqsec1-alekandrov-variation-formula}
\left.\frac{d}{dt}\right|_{t=0}V([h_0+tf])
=
\int_{\Sn}f(v)\,dS([h_0],v);
\end{equation}
see \cite{Aleksandrov1939}.Formula \eqref{eqsec1-alekandrov-variation-formula}
reduces the Minkowski problem to an optimization problem
and identifies the surface area measure as the differential
of the volume functional. A key feature of this formula is that the derivative is two-sided and the formula remains valid without smoothness or strict convexity assumptions on the convex body.

The same variational approach applies to other geometric
invariants. Jerison \cite{Jerison1996} studied the capacitary
measure arising from the first variation of electrostatic capacity and established existence and regularity
results for the corresponding Minkowski problem. Huang--Lutwak--Yang--Zhang \cite{Huang2016Acta} proved that variations of the dual quermassintegrals produce the dual curvature measures, thereby identifying the geometric measures underlying the dual Brunn--Minkowski theory and formulating the associated dual Minkowski problem. These examples illustrate the principle that variational formulas for natural geometric invariants give rise to geometric measures and provide a natural framework for studying the associated Minkowski problems.

 However, not every natural geometric invariant admits a
 comparably robust variational formula. Let $K,L\subset\Rn$
 be convex bodies, and let
 $K+L=\{x+y:x\in K,\ y\in L\}$ denote their Minkowski addition.
 Write $S(K)=\cHnmone(\partial K)$ for the surface area of $K$.
 The first variation of surface area under Minkowski addition is
 \begin{equation}\label{eqsec1-variation-of-surface-area}
 	\left.\frac{d}{dt}\right|_{t=0^+}S(K+tL)
 	=
 	(n-1)\int_{\Sn}h_L(v)\,dS_{n-2}(K,v),
 \end{equation}
 where $h_L$ is the support function of $L$ and
 $S_{n-2}(K,\cdot)$ is the $(n-2)$-th area measure of $K$;
 see \cite{schneiderbook}. Unlike \eqref{eqsec1-alekandrov-variation-formula}, the derivative in \eqref{eqsec1-variation-of-surface-area} is one-sided and is formulated along Minkowski additions of convex bodies. The Minkowski problems for area measures are called the Christoffel--Minkowski problems. For the solution of the Christoffel problem,
 see Firey \cite{Firey1968}, Berg \cite{Berg1969}, Grinberg-Zhang \cite{GrinbergZhang1999} and Goodey-Yaskin-Yaskina \cite{GoodeyYaskinYaskina2011}.
For contributions to the Christoffel--Minkowski problem
and related problems using PDE methods,
see Guan--Guan \cite{GuanGuan2002},
Guan--Li--Li \cite{GuanLiLi2012},
Guan--Ma \cite{GuanMa2003}, and
Guan--Ma--Zhou \cite{GuanMaZhou2006}.
However, the general Christoffel--Minkowski
problems remain open.
 
Recently, Lutwak--Xi--Yang--Zhang \cite{XiCPAM} used
variational methods to study the chord integrals in integral geometry. Let $\mathscr{L}^n$ denote the affine Grassmannian of lines in $\Rn$, equipped with the normalized invariant measure. For $q\ge0$, the $q$-th chord integral of a convex body $K$ is
\[
I_q(K)=\int_{\mathscr{L}^n}|K\cap\ell|^q\,d\ell,
\]
where $|K\cap\ell|$ denotes the length of the chord $K\cap\ell$. Using Crofton's volume formula, Cauchy's formula for surface area, and the Poincar\'e-Hadwiger formula, respectively (see \cite{Ren1994,Santalo2004}), we have
\[
I_1(K)=V(K),\qquad
I_0(K)=\frac{\omega_{n-1}}{n\omega_n}S(K),\qquad
I_{n+1}(K)=\frac{n+1}{\omega_n}V(K)^2.
\]
Thus, the chord integrals contain volume and surface area as special cases. For $q>0$, Lutwak--Xi--Yang--Zhang established robust
Aleksandrov-type variational formulas for $I_q$, giving
rise to the chord measures $F_q(K,\cdot)$. They formulated and solved the corresponding Minkowski problems. Moreover, under suitable regularity assumptions on $K$,
the limiting measure of $F_q(K,\cdot)$ as $q\to0^+$ is
\[
\frac{(n-1)\omega_{n-1}}{n\omega_n}S_{n-2}(K,\cdot).
\]
Chord integrals of convex bodies can also be expressed through the $X$-ray functions. For $u\in\Sn$ and $y\in u^\perp$, define
\[
X_K(y,u)=\cH^1\bigl(K\cap(y+\mathbb Ru)\bigr).
\]
Then, for $q>0$,
\[
I_q(K)=\frac{1}{n\omega_n}\int_{\Sn}\int_{u^\perp}X_K(y,u)^q\,dy\,du.
\]
For $q=0$, the inner integral is replaced by $\vol_{n-1}(K|u^\perp)$; see \cite{XiCPAM}.

%

Motivated by this representation, we introduce the $(s,q)$-chord integral. For $q>0$ and $s\in\bR\setminus\{0\}$, define
\[
\Psi_{s,q}(K)
=
\frac{1}{n\omega_n}
\int_{\Sn}
\left(\int_{u^\perp}X_K(x,u)^q\,dx\right)^s\,du.
\]
When $q=0$, set
\[
\Psi_{s,0}(K)
=
\frac{1}{n\omega_n}\int_{\Sn}\vol_{n-1}(K|u^\perp)^s\,du.
\]
The complete definition is given in Definition~\ref{Def-sq-integral}. When $s=1$,  $\Psi_{1,q}(K)$ is the chord integral $I_q(K)$ of $K$. 

We also introduce a new geometric measures called the
$(s,q)$-chord measure $F_{s,q}(K,\cdot)$ for a convex body $K$;
see Definition~\ref{Def-sq-chord-measure}. The following list shows that the $(s,q)$-chord measures cover many important geometric measures as special cases.
\begin{itemize}
	\item When $q=1$, we have
	\[
	F_{s,1}(K,\cdot)=V(K)^{s-1}S(K,\cdot).
	\]
	In particular, $F_{1,1}(K,\cdot)$ is the classical surface
	area measure.
	
	\item When $s=1$, the $(s,q)$-chord measures recover the chord measures
	\[
	F_{1,q}(K,\cdot)=F_q(K,\cdot),
	\qquad q>0.
	\]
	
	\item When $s=n/(q-1)$, the $(s,q)$-chord measures recover, up to normalization,
	the fractional affine area measures introduced by
	Xi--Zhao \cite{XiZhao26}:
	\[
	F_{\frac{n}{q-1},q}(K,\cdot)
	=\frac{2}{\omega_n}F^{(a)}_{q-1}(K,\cdot),
	\qquad q>0,\quad q\ne1.
	\]
\end{itemize}
  Our first result shows that, for $q>0$, the variation
  of $\Psi_{s,q}(K)$ leads to the $(s,q)$-chord measure
  $F_{s,q}(K,\cdot)$.
 For a compact set $\Omega$, let $C(\Omega)$ be the set of  real-valued continuous functions on $\Omega$ and let $C^+(\Omega)$ be the set of positive functions in $C(\Omega)$. Let $o(t,\cdot)\in C(\Omega)$ denote a continuous function satisfying $o(t,\cdot)/t\to0$ uniformly as $t\to0$.

\begin{theorem}\label{thm-variation-formula-for-sq-chord-int}
Let $q>0$, $s\in\bR\setminus\{0\}$, and let $\Omega\subset\Sn$ be compact and not contained in any closed hemisphere. Let $h_0,g\in C(\Omega)$. Suppose that, for sufficiently small $\delta>0$ and $|t|<\delta$, the Wulff shape $K_t=[h_t]$ is a convex body in $\Rn$, where
\[
h_t=h_0+tg+o(t,\cdot)\qquad\text{on }\Omega.
\]
Then
\[
\left.\frac{d}{dt}\right|_{t=0}\Psi_{s,q}(K_t)
=s\int_\Omega g(v)\,dF_{s,q}(K_0,v).
\]
\end{theorem}

With this variation formula, we study the Minkowski problems associated with $F_{s,q}(K,\cdot)$ and obtain the following result.

\begin{theorem}\label{thm-mp-for-sq-chord-measure}
Let $q>0$ and $s\in\bR\setminus\{0\}$ satisfy $(n+q-1)s\ne1$. If $\mu$ is a nonzero finite Borel measure on $\Sn$, then there exists a convex body $K$ such that
\[
F_{s,q}(K,\cdot)=\mu
\]
if and only if $\mu$ is not concentrated on any closed hemisphere of $\Sn$ and
\begin{equation}\label{eq-centroid-of-mu}
\int_{\Sn}v\,d\mu(v)=o.
\end{equation}
\end{theorem}
When $s<0$, we establish a Brunn--Minkowski-type inequality
for the $(s,q)$-chord integrals; see Theorem~\ref{thm:secondBM}.
Using this inequality, we show that the solution is unique
up to translation; see Corollary~\ref{coro-unique-of-mp}.

An important extension of the classical Minkowski problem is the $L_p$ Minkowski problem, introduced by Lutwak \cite{Lutwak1993}. The case $p=1$ is the classical Minkowski problem. The existence, uniqueness, regularity, and discrete versions for $p>1$ were developed in a series of works; see, for example, \cite{LutwakOliker1995,LutwakYangZhang2004,HugLutwakYangZhang2005,ChouWang2006}. The range $p<1$ is more delicate, since uniqueness may fail and the behavior of solutions depends strongly on $p$ and on the prescribed measure. Various existence, nonuniqueness, and regularity results in this range were established in \cite{JianLuWang2015,BoroczkyTrinh2017,ChenLiZhu2017,Zhu2017,BianchiBoroczkyColesantiYang2019,Chen2006,BianchiBoroczkyColesanti2020}. In particular, when $p=0$, the corresponding problem is the logarithmic Minkowski problem. It asks for necessary and sufficient conditions on a finite
Borel measure $\mu$ on $\Sn$ such that
\[
d\mu(v)=h_K(v)\,dS(K,v)
\]
for some $K\in\mathcal K^n_o$.
The logarithmic Minkowski problem is closely connected with the log-Brunn--Minkowski
inequality \cite{BoroczkyLutwakYangZhang2012}. See \cite{Boroczky2013JAMS,Zhu2014,BoroczkyHegedusZhu2016,ChenLiZhu2019,ChenFengLiu2022} for works on the logarithmic Minkowski problem.

Motivated by the above, we introduce the $L_p$ extension of
the $(s,q)$-chord measure (see Definition~\ref{def-Lp-sq-chord-measure})
and the associated $L_p$ Minkowski problem. In the smooth case, the corresponding measure equations reduce to a nonlocal Monge--Amp\`ere equation \eqref{eqsec3-MA-equation}. The cases $p=1$ and $p=0$ correspond, respectively, to the ordinary and logarithmic Minkowski problems for $(s,q)$-chord measures. We study the latter in the even case and obtain the following result.

\begin{theorem}\label{thmsec4-solution-to-log-mp}
Let $q\ge1$ and $s>0$ satisfy $(q-1)s\in[0,n]$, and let $\mu$ be a nonzero even finite Borel measure on $\Sn$. Suppose that
\[
\frac{\mu(\xi\cap\Sn)}{|\mu|}
<
\frac{ks+\min\{k,(q-1)s\}}{(n+q-1)s}
\]
for every $k\in\{1,\ldots,n-1\}$ and every $k$-dimensional linear subspace $\xi\subset\Rn$. Then there exists an origin-symmetric convex body $K$ such that
\[
d\mu(v)=h_K(v)\,dF_{s,q}(K,v).
\]
\end{theorem}

Theorem~\ref{thmsec4-solution-to-log-mp} extends the result
on the logarithmic Minkowski problem for chord measures
in \cite{XiCPAM} to $(s,q)$-chord measures. A key ingredient in its proof is the $(s,q)$-chord estimate for ellipsoids established in Lemma~\ref{lemsec4-ellipsoid-estimate}. The parameter $s$ in the $(s,q)$-chord integral causes
substantial difficulties for the estimate.
We use the star body $S_{q-1}(K)$,
which was introduced by Gardner--Zhang \cite{GardnerAJM98}, and the affine covariance of $S_{q-1}(K)$
(see \eqref{eqsec4-def-of-S-q-minus-one-K} and
\eqref{eqsec4-affine-property-of-S-q-minus-one-K}),
to overcome these difficulties.

This paper is organized as follows. In Section~\ref{sec:preliminaries}, we introduce the necessary preliminaries. In Section~\ref{sec:variation}, we prove Theorem~\ref{thm-variation-formula-for-sq-chord-int} and establish finiteness, translation invariance, homogeneity, and the centroid property of the $(s,q)$-chord measure. In Section~\ref{sec:minkowski}, we prove Theorem~\ref{thm-mp-for-sq-chord-measure}. In Section~\ref{sec:log-minkowski}, we prove Theorem~\ref{thmsec4-solution-to-log-mp}. Finally, in Section~\ref{sec:asymptotic}, we study the asymptotic behavior of $F_{s,q}(K,\cdot)$ as $q\to0^+$ under suitable regularity assumptions on $\partial K$.

\section{Preliminaries}\label{sec:preliminaries}

We continue to use the conventions introduced above. 
The orthogonal projection of $K\subset\Rn$ onto $u^\perp$ is denoted by $K|u^\perp$. If $A\subset\mathbb R^n$ is contained in an $m$-dimensional affine subspace, then $\vol_m(A)$ denotes its $m$-dimensional Lebesgue measure. The Grassmannian of $k$-dimensional linear subspaces of $\Rn$ is denoted by $G_{n,k}$.
A finite Borel measure $\mu$ on $\Sn$ has centroid at the origin if
\[
\int_{\Sn}v\,d\mu(v)=o.
\]
It is not concentrated on any closed hemisphere if
\[
\mu(\{v\in\Sn:v\cdot u\ge0\})<|\mu|
\qquad\text{for every }u\in\Sn.
\]
The centroid of $K\in\Kn$ is defined by
\[
c(K)=\frac{1}{V(K)}\int_K x\,dx.
\]
For $A,B\subset\Rn$ and $t>0$, let
\[
A+B=\{a+b:a\in A,\ b\in B\},
\]
and\[
tB=\{tb:b\in B\}.\]
We write $A+x=A+\{x\}$ for $x\in\mathbb R^n$. Convex bodies $K,L$ are homothetic if $K=tL+x$ for some $t>0$ and $x\in\Rn$; they are translates if $t=1$.

For a compact convex set $K$, its support function $h_K:\Sn\to\bR$ is
\begin{equation*}
h_K(v)=\max\{x\cdot v:x\in K\},
\qquad v\in\Sn.
\end{equation*}
Its width function in direction $u\in\Sn$ is
\begin{equation*}
	w_K(u)=h_K(u)+h_K(-u).
\end{equation*}
For compact convex sets $K_i,K\subset\Rn$, the sequence $(K_i)$
converges to $K$ in the Hausdorff metric if
\[
\|h_{K_i}-h_K\|_\infty
=\max_{u\in\Sn}|h_{K_i}(u)-h_K(u)|
\longrightarrow0.
\]
The projection body $\Pi K$ of $K\in\Kn$ is the origin-symmetric convex body with support function
\[
h_{\Pi K}(u)=\vol_{n-1}(K|u^\perp),
\qquad u\in\Sn.
\]
For $K\in\Kn$ and $x\in\partial K$, a vector $v\in\Sn$ is called an outer unit normal of $\partial K$ at $x$ if
\[
x\cdot v=h_K(v).
\]
Let $\partial'K$ denote the set of boundary points of $K$ that  have a unique outer unit normal.
It is well known \cite[Theorem~2.2.5]{schneiderbook} that \[\cH^{n-1}(\partial K\setminus\partial'K)=0.\]
 For $x\in\partial'K$, the Gauss map $\nu_K(x)$ of $K$ is defined as the unique outer unit normal of $\partial K$ at $x$. 
The surface area measure $S(K,\cdot)$ is defined by
\begin{equation*}
	S(K,\omega)=\cHnmone(\nu_K^{-1}(\omega)),
\end{equation*}where $\omega\subset\Sn$ is an arbitrary Borel set.

Let $\Omega\subset\Sn$ be compact and not contained in any closed hemisphere. For $f\in C(\Omega)$, the Wulff shape generated by $f$ is
\begin{equation}\label{defsec2-wulff-shape}
	[f]=\bigcap_{v\in\Omega}\{x\in\Rn:v\cdot x\le f(v)\}.
\end{equation}
Whenever $[f]$ is nonempty, $h_{[f]}\le f$ on $\Omega$. If $K=[f]\in\Kn$, then
\(
h_K(v)=f(v)
\) for $S(K,\cdot)$-almost every $v\in\Omega$. If $K\in\Kn$ and $\Omega=\Sn$, then $[h_K]=K$.

A star-shaped set $Q$ with respect to the origin is a compact
set whose intersection with every ray emanating from $o$
is a line segment with one endpoint at $o$.
If $K$ is a star-shaped set with respect to the origin,
its radial function is
\[
\rho_K(u)=\max\{t\ge0:tu\in K\},
\qquad u\in\Sn.
\]
 A star-shaped set is called a star body if its radial function is positive and continuous on $\Sn$. 
 More generally, if $y\in K$ and $x\in\Rn\setminus\{o\}$, the extended radial function of $K$ with respect to $y$ is
\[
\rho_{K,y}(x)=\max\{t\ge0:y+tx\in K\}.
\]

We recall that the $X$-ray of a convex body $K$ is
\[
X_K(y,u)=\cH^1\bigl(K\cap(y+\mathbb Ru)\bigr),
\qquad y\in\Rn,\quad u\in\Sn.
\]
If $z\in K\cap(y+\mathbb Ru)$, then
\begin{equation}\label{eq-Xrho}
X_K(y,u)=\rho_{K,z}(u)+\rho_{K,z}(-u).
\end{equation}
 The directional chord integral $I_q(K,u)$ for $q>0$ and $u\in\Sn$ is defined by 
\begin{equation*}
I_q(K,u)
=
\int_{u^\perp}X_K(x,u)^q\,dx.
\end{equation*}
For $q=0$, set
\begin{equation*}
I_0(K,u)=\vol_{n-1}(K|u^\perp).
\end{equation*}
For $x_0\in\mathbb R^n$, let $x_0|u^\perp$ be the orthogonal projection of $x_0$ on $u^\perp$. Since \[X_{K+x_0}(y,u)=X_{K}(y-x_0|u^\perp,u),\]by a change of variables, we have \begin{equation}\label{eqsec2 translation invariance of I_q(K,u)}
	I_q(K,u)=I_q(K+x_0,u).
\end{equation} Let $\lambda>0$. Then we have
\begin{equation}\label{eq-homogeneity-of-Iq}
	I_q(\lambda K,u)=\lambda^{n+q-1}I_q(K,u).
\end{equation}
 The function $u\mapsto I_q(K,u)$ is positive and continuous on $\Sn$ for $q\ge0$; for $q>0$, see \cite{GardnerAJM98}.
 For $K\in\Kn$ and $q\in[0,\infty)\setminus\{1\}$, Gardner-Zhang \cite{GardnerAJM98} introduced the star body $S_{q-1}K$ which is defined by 
 \begin{equation}\label{eqsec4-def-of-S-q-minus-one-K}
 	\rho_{S_{q-1}K}(u)=I_q(K,u)^{1/(q-1)},
 	\qquad u\in\Sn.
 \end{equation}
 The star body $S_{q-1}K$ is affine covariant. For $\phi\in\operatorname{SL}(n)$, $x_0\in\Rn$, and $\lambda>0$,
 \begin{equation}\label{eqsec4-affine-property-of-S-q-minus-one-K}
 	S_{q-1}(\phi K+x_0)=\phi S_{q-1}K,
 	\qquad
 	S_{q-1}(\lambda K)=\lambda^{(n+q-1)/(q-1)}S_{q-1}K,
 \end{equation}
 see \cite{XiZhao26}.
 
Now we define the $(s,q)$-chord integral. 
\begin{definition}\label{Def-sq-integral}
Let $K\in\Kn$, $q\ge0$, and $s\in\bR\setminus\{0\}$. The $(s,q)$-chord integral of $K$ is
\begin{equation}\label{def-sq-chord-int}
\Psi_{s,q}(K)
=
\frac{1}{n\omega_n}\int_{\Sn}I_q(K,u)^s\,du.
\end{equation}
\end{definition}
Since $K\in\mathcal K^n$, $\Psi_{s,q}(K)$ is positive and finite. By \eqref{eqsec2 translation invariance of I_q(K,u)}, we know that $\Psi_{s,q}(K)$ is translation invariant. The $(s,q)$-chord integral is also continuous in the Hausdorff metric. Indeed, suppose that $K_i\to K$ in the Hausdorff metric. Since the $(s,q)$-chord integral is translation invariant, after translating all bodies by a fixed point, we can assume that $K\in\mathcal K_o^n$. So there are numbers $\varepsilon_i\to0^+$ such that, for all sufficiently large $i$,
\[
(1-\varepsilon_i)K\subset K_i\subset(1+\varepsilon_i)K.
\]
Monotonicity under inclusion and homogeneity of the directional chord integrals then give, uniformly in $u\in\Sn$,
\[
(1-\varepsilon_i)^{n+q-1}I_q(K,u)
\le I_q(K_i,u)
\le(1+\varepsilon_i)^{n+q-1}I_q(K,u).
\] From this fact, the continuity of $\Psi_{s,q}$ is immediate.
$\Psi_{s,q}$ also has the following homogeneity. 

\begin{lemma}\label{lem-homogeneity-of-sq-chord-int}
Let $\lambda>0$, $q\ge0$, $s\in\bR\setminus\{0\}$, and $K\in\Kn$. Then
\[
\Psi_{s,q}(\lambda K)=\lambda^{(n+q-1)s}\Psi_{s,q}(K).
\]
\end{lemma}

\begin{proof}
The assertion follows immediately from \eqref{eq-homogeneity-of-Iq} and Definition~\ref{Def-sq-integral}.
\end{proof}

 Now we define the $(s,q)$-chord measures.
\begin{definition}\label{Def-sq-chord-measure}
Let $q>0$, $s\in\bR\setminus\{0\}$, and $K\in\Kn$. The $(s,q)$-chord measure of $K$ is the Borel measure on $\Sn$ defined by
\begin{equation}\label{def-sq-chord-measure}
F_{s,q}(K,\omega)
=
\frac{q}{n\omega_n}
\int_{\nu_K^{-1}(\omega)}\int_{\Sn}
I_q(K,u)^{s-1}X_K(y,u)^{q-1}\,du\,d\cHnmone(y)
\end{equation}
for every Borel set $\omega\subset\Sn$.
\end{definition}
 \begin{remark}
 	Let
 	\[
 	A_K=\{(y,u)\in\partial'K\times\Sn:X_K(y,u)=0\}.
 	\]
 	For each $y\in\partial'K$, 
 	$A_K(y)=\{u\in\Sn:(y,u)\in A_K\}$ is contained in
 	$\Sn\cap\nu_K(y)^\perp$, which has zero spherical measure.
 	Thus, for every Borel set $\omega\subset\Sn$,
 	\[
 	\int_{\nu_K^{-1}(\omega)}\int_{\Sn}
 	\mathbf{1}_{A_K}(y,u)\,du\,d\cHnmone(y)=0.
 	\]
 	Consequently, the integral in \eqref{def-sq-chord-measure}
 	can be taken over the pairs $(y,u)$ for which $X_K(y,u)>0$,
 	without changing its value.
 \end{remark}


Let $K\in\Kn$ and $z\in K$. For $q>-1$, the $q$-th dual quermassintegral is defined, with values in $[0,+\infty]$, by
\[
\widetilde V_q(K,z)
=
\frac{1}{n\omega_n}
\int_{S^+_z(K)}
\rho_{K,z}(u)^q\,du,
\]
where
\[
S^+_z(K)=\{u\in\Sn:\rho_{K,z}(u)>0\}.
\]
When $o\in K$, we write
\[
\widetilde V_q(K)=\widetilde V_q(K,o).
\]If $q>0$, then
\begin{equation}\label{eqsec5-dual-quermassintegral-by-integration-on-K}
	\widetilde V_q(K,z)
	=
	\frac{q}{n\omega_n}\int_K |x-z|^{q-n}\,dx.
\end{equation}
This identity follows by using polar coordinates. For later convenience, our definition includes an additional factor
of $1/\omega_n$ compared with the usual normalization of the
dual quermassintegral.

For $q>0$, the chord measures $F_q(K,\cdot)$ of Lutwak--Xi--Yang--Zhang in \cite{XiCPAM} are defined by
\begin{equation}\label{def-classical-chord-measure}
	F_q(K,\omega)
	=
	\frac{2q}{n\omega_n}
	\int_{\nu_K^{-1}(\omega)}
	\int_{S^+_z(K)}\rho_{K,z}(u)^{q-1}\,du\,d\cHnmone(z)
\end{equation}
for every Borel set $\omega\subset\Sn$.

When $s=1$, the $(s,q)$-chord measure recovers the chord measure, namely
\begin{equation*}
	F_{1,q}(K,\cdot)=F_q(K,\cdot).
\end{equation*}
Indeed, for each $z\in\partial'K$, the sets $S_z^+(K)$ and
$-S_z^+(K)$ form a partition of $\Sn$ up to a set of spherical
measure zero. Thus, by \eqref{eq-Xrho} and the change of variables
$u\mapsto -u$,
\begin{align*}
	\int_{\Sn}X_K(z,u)^{q-1}\,du
	&=
	\int_{S_z^+(K)}\rho_{K,z}(u)^{q-1}\,du
	+
	\int_{-S_z^+(K)}\rho_{K,z}(-u)^{q-1}\,du\\
	&=
	2\int_{S_z^+(K)}\rho_{K,z}(u)^{q-1}\,du.
\end{align*}
Substituting this identity into \eqref{def-sq-chord-measure}
with $s=1$ gives \eqref{def-classical-chord-measure}.

\section{\texorpdfstring{The variational formula for the $(s,q)$-chord integral}{The variational formula for the (s,q)-chord integral}}\label{sec:variation}

In this section we prove Theorem~\ref{thm-variation-formula-for-sq-chord-int} and record several basic properties of $F_{s,q}(K,\cdot)$. The first point is to show that $F_{s,q}(K,\cdot)$ is indeed a finite Borel measure for every convex body $K$. We need the following theorem.

\begin{theorem}[Theorem 3.3, \cite{XiZhao26}]\label{prop-boundary-representation-of-Iq}
Let $K\in\Kn$ and $q>0$. If $\partial K$ is line-free in $u\in\Sn$, then
\[
I_q(K,u)
=
\frac{q}{n+q-1}
\int_{\partial K}X_K(z,u)^{q-1}\bigl(z\cdot\nu_K(z)\bigr)\,d\cHnmone(z).
\]
\end{theorem}
We say that $\partial K$ is line-free in the direction $u\in\Sn$ if it contains no nontrivial line segment parallel to $u$. For a convex body, $\partial K$ is line-free in almost every direction.
\begin{lemma}\label{lemsec2-finite-of-Fsq}
Let $K\in\Kn$, $q>0$, and $s\in\bR\setminus\{0\}$. Then
\[
F_{s,q}(K,\Sn)<\infty.
\]
\end{lemma}

\begin{proof} A direct change of variables in Definition~\ref{Def-sq-chord-measure} shows that $F_{s,q}(K,\Sn)$ is invariant under translations of $K$. We may therefore assume that $K\in\Kno$. So there exists a constant $c_0>0$
	 such that \begin{equation}\label{eqsec2-positive-support-bound} y\cdot\nu_K(y)\ge c_0 \end{equation} for $\cHnmone$-almost every $y\in\partial K$. 
	 Since the function $u\longmapsto I_q(K,u)^{s-1}$ is continuous and positive, there exists a constant $C>0$ such that \begin{equation}\label{eq:Iq-weight-bound} I_q(K,u)^{s-1}\le C \qquad\text{for every }u\in\Sn. \end{equation} 
	    By Definition~\ref{Def-sq-chord-measure} and \eqref{eq:Iq-weight-bound}, we have \begin{align} F_{s,q}(K,\Sn) &= \frac{q}{n\omega_n} \int_{\partial K}\int_{\Sn} I_q(K,u)^{s-1} X_K(y,u)^{q-1} \,du\,d\cHnmone(y) \notag\\ &\le \frac{qC}{n\omega_n} \int_{\partial K}\int_{\Sn} X_K(y,u)^{q-1} \,du\,d\cHnmone(y). \label{eq:Fsq-first-bound} \end{align} By \eqref{eqsec2-positive-support-bound}, \[ 1 \le \frac{y\cdot\nu_K(y)}{c_0} \] for $\cHnmone$-almost every $y\in\partial K$. Hence \begin{align} F_{s,q}(K,\Sn) &\le \frac{qC}{n\omega_nc_0} \int_{\partial K}\int_{\Sn} X_K(y,u)^{q-1} \bigl(y\cdot\nu_K(y)\bigr) \,du\,d\cHnmone(y) \notag\\ &= \frac{qC}{n\omega_nc_0} \int_{\Sn} \int_{\partial K} X_K(y,u)^{q-1} \bigl(y\cdot\nu_K(y)\bigr) \,d\cHnmone(y) du, \label{eq:Fsq-tonelli} 
	    \end{align} 
	    where the last equality follows from Tonelli's theorem. For almost every $u\in\Sn$, Theorem~\ref{prop-boundary-representation-of-Iq} applies.  We obtain \begin{align*} F_{s,q}(K,\Sn) \le \frac{(n+q-1)C}{n\omega_nc_0} \int_{\Sn}I_q(K,u)\,du. \end{align*} Finally, $I_q(K,\cdot)$ is continuous on $\Sn$. Consequently, \[ \int_{\Sn}I_q(K,u)\,du  <\infty. \] 
	    It follows that \[ F_{s,q}(K,\Sn)<\infty, \] as claimed. \end{proof}

In order to prove Theorem \ref{thm-variation-formula-for-sq-chord-int}, we shall use the following robust variation formula. 

\begin{theorem}[Theorem 3.7, \cite{XiZhao26}]\label{prop-variation-formula-for-Iq}
Let $q>0$ and let $K_t$ be as in Theorem~\ref{thm-variation-formula-for-sq-chord-int}. Then, for almost every $u\in\Sn$,
\[
\left.\frac{d}{dt}\right|_{t=0}I_q(K_t,u)
=
q\int_{\partial K_0}
g(\nu_{K_0}(z))X_{K_0}(z,u)^{q-1}\,d\cHnmone(z).
\]
\end{theorem}

\begin{proof}[Proof of Theorem~\ref{thm-variation-formula-for-sq-chord-int}]
The functional $\Psi_{s,q}$ and every directional chord integral are translation invariant. After translating, we may assume that $K_0\in\Kno$; equivalently, $h_0\in C^+(\Omega)$. Since $g$ is continuous, there is a constant $c>0$ such that, for all sufficiently small $|t|$,
\begin{equation*}
(1-c|t|)K_0\subset K_t\subset(1+c|t|)K_0.
\end{equation*}
By monotonicity of $I_q$ under inclusion and \eqref{eq-homogeneity-of-Iq}, for either sign of $s$,
\begin{align*}
\left|
\frac{I_q(K_t,u)^s-I_q(K_0,u)^s}{t}
\right|
&\le
\left|
\frac{I_q((1+c|t|)K_0,u)^s-I_q((1-c|t|)K_0,u)^s}{t}
\right|\\
&=
\left|
\frac{(1+c|t|)^{(n+q-1)s}-(1-c|t|)^{(n+q-1)s}}{t}
\right|I_q(K_0,u)^s.
\end{align*}
Note that there exists $C_0>0$ such that for sufficiently small $|t|\neq0$,\[\left|
\frac{(1+c|t|)^{(n+q-1)s}-(1-c|t|)^{(n+q-1)s}}{t}
\right|\le C_0.\] Since $I_q(K_0,\cdot)^s$ is integrable, dominated convergence and Theorem~\ref{prop-variation-formula-for-Iq} therefore give
\begin{align*}
\left.\frac{d}{dt}\right|_{t=0}\Psi_{s,q}(K_t)
&=
\frac{sq}{n\omega_n}
\int_{\Sn}I_q(K_0,u)^{s-1}
\int_{\partial K_0}g(\nu_{K_0}(z))X_{K_0}(z,u)^{q-1}\,d\cHnmone(z)\,du\\
&=
\frac{sq}{n\omega_n}
\int_{\partial K_0}g(\nu_{K_0}(z))
\int_{\Sn}I_q(K_0,u)^{s-1}X_{K_0}(z,u)^{q-1}\,du\,d\cHnmone(z)\\
&=
s\int_\Omega g(v)\,dF_{s,q}(K_0,v).
\end{align*}
Here the last equality uses Definition~\ref{Def-sq-chord-measure} and the fact that the surface area measure of the Wulff shape $K_0=[h_0]$ is concentrated on $\Omega$.
\end{proof}

The following theorem shows that $\Psi_{s,q}(K)$ can be
recovered from $F_{s,q}(K,\cdot)$.

\begin{theorem}\label{thm-integral-representation-of-sq-chord}
Let $K\in\Kn$, $q>0$, and $s\in\bR\setminus\{0\}$. Then
\[
\int_{\Sn}h_K(v)\,dF_{s,q}(K,v)
=(n+q-1)\Psi_{s,q}(K).
\]
\end{theorem}
\begin{proof}
	We apply Theorem~\ref{thm-variation-formula-for-sq-chord-int}
	with
	\[
	\Omega=\Sn,\qquad h_0=h_K,\qquad g=h_K,
	\]
	and
	\[
	h_t=(1+t)h_K.
	\]
	Since
	\[
	h_t=h_K+t h_K,
	\]
	Theorem~\ref{thm-variation-formula-for-sq-chord-int} gives
	\begin{equation}\label{eq-variation-for-dilation}
		\left.\frac{d}{dt}\right|_{t=0}
		\Psi_{s,q}((1+t)K)
		=
		s\int_{\Sn}h_K(v)\,dF_{s,q}(K,v).
	\end{equation}
	On the other hand, by
	Lemma~\ref{lem-homogeneity-of-sq-chord-int},
	\[
	\Psi_{s,q}((1+t)K)
	=
	(1+t)^{(n+q-1)s}\Psi_{s,q}(K).
	\]
	Differentiating at $t=0$, we obtain
	\begin{align}
		\left.\frac{d}{dt}\right|_{t=0}
		\Psi_{s,q}((1+t)K)=
		(n+q-1)s\,\Psi_{s,q}(K).
		\label{eq-homogeneity-derivative}
	\end{align}
	Comparing \eqref{eq-variation-for-dilation} and
	\eqref{eq-homogeneity-derivative}, and using $s\neq0$, we obtain
	\[
	\int_{\Sn}h_K(v)\,dF_{s,q}(K,v)
	=
	(n+q-1)\,\Psi_{s,q}(K).
	\]
	This proves the theorem.
\end{proof}
The following proposition shows that the \((s,q)\)-chord measure is translation invariant.
\begin{proposition}\label{pro-translation-invariant}
Let $K\in\Kn$, $q>0$, and $s\in\bR\setminus\{0\}$. Then
\[
F_{s,q}(K+x_0,\cdot)=F_{s,q}(K,\cdot)
\qquad\text{for every }x_0\in\Rn.
\]
\end{proposition}

\begin{proof}
	Let $\Omega=\Sn$ and $h_0=h_K$ in
	Theorem~\ref{thm-variation-formula-for-sq-chord-int}. Then, for every
	$g\in C(\Sn)$,
	\[
	\left.\frac{d}{dt}\right|_{t=0}\Psi_{s,q}(K_t)
	=
	s\int_{\Sn}g(v)\,dF_{s,q}(K,v).
	\]
Taking $h_0=h_{K+x_0}$, we obtain
	\[
	\left.\frac{d}{dt}\right|_{t=0}\Psi_{s,q}(K_t+x_0)
	=
	s\int_{\Sn}g(v)\,dF_{s,q}(K+x_0,v).
	\]
	Since $\Psi_{s,q}$ is translation invariant,
	\[
	\Psi_{s,q}(K_t+x_0)=\Psi_{s,q}(K_t).
	\]
	Thus the above two identities, together with $s\neq0$, imply
	\[
	\int_{\Sn}g(v)\,dF_{s,q}(K+x_0,v)
	=
	\int_{\Sn}g(v)\,dF_{s,q}(K,v)
	\]
	for every $g\in C(\Sn)$. Hence the assertion follows.
\end{proof}

We now show that the centroid of the $(s,q)$-chord measure
is at the origin.
\begin{theorem}\label{thm-centroid-of-sq-chord-measure}
Let $K\in\Kn$, $q>0$, and $s\in\bR\setminus\{0\}$. Then
\[
\int_{\Sn}v\,dF_{s,q}(K,v)=o.
\]
\end{theorem}

\begin{proof}
Let $x_0\in\Rn$. By Theorem~\ref{thm-integral-representation-of-sq-chord}, translation invariance of $\Psi_{s,q}$, and Proposition~\ref{pro-translation-invariant},
\begin{align*}
\int_{\Sn}h_K(v)\,dF_{s,q}(K,v)
&=
\int_{\Sn}h_{K+x_0}(v)\,dF_{s,q}(K+x_0,v)\\
&=
\int_{\Sn}\bigl(h_K(v)+x_0\cdot v\bigr)\,dF_{s,q}(K,v).
\end{align*}
Hence
\[
x_0\cdot\int_{\Sn}v\,dF_{s,q}(K,v)=0.
\]
Since $x_0$ is arbitrary, the centroid of $F_{s,q}(K,\cdot)$ lies at the origin.
\end{proof}

Now we introduce the $L_p$ $(s,q)$-chord measures.
\begin{definition}\label{def-Lp-sq-chord-measure}
	Let $K\in\Kno$, $q>0$, $s\in\bR\setminus\{0\}$, and $p\in\bR$. The $L_p$ $(s,q)$-chord measure is defined by
	\[
	dF_{p,s,q}(K,v)=h_K(v)^{1-p}\,dF_{s,q}(K,v),
	\qquad v\in\Sn.
	\]
\end{definition}
By Theorem~\ref{thm-variation-formula-for-sq-chord-int},
we obtain the following variational formula for the
$L_p$ $(s,q)$-chord measure.
\begin{theorem}\label{thmsec2-variational-formula-for-Lp}
	Let $q>0$, $p\in\bR$, $s\in\bR\setminus\{0\}$, and let $\Omega\subset\Sn$ be compact and not contained in any closed hemisphere. Let $h_0\in C^+(\Omega)$ and $g\in C(\Omega)$. Suppose that, for sufficiently small $|t|$, $K_t=[h_t]\in\Kno$ and
	\[
	\begin{cases}
		h_t^p=h_0^p+tg+o(t,\cdot),&p\ne0,\\
		\log h_t=\log h_0+tg+o(t,\cdot),&p=0.
	\end{cases}
	\]
	Then
	\[
	\left.\frac{d}{dt}\right|_{t=0}\Psi_{s,q}(K_t)
	=
	\begin{cases}
		\displaystyle \frac{s}{p}\int_\Omega g(v)\,dF_{p,s,q}(K_0,v),&p\ne0,\\[6pt]
		\displaystyle s\int_\Omega g(v)\,dF_{0,s,q}(K_0,v),&p=0.
	\end{cases}
	\]
\end{theorem}


We have the following homogeneity for the $L_p$ $(s,q)$-chord measures.
\begin{proposition}\label{prop-homogeneity-of-sq-chord-measure}
	Let $\lambda>0$, $q>0$, $s\in\bR\setminus\{0\}$, and $K\in\Kn$. Then
	\[
	F_{s,q}(\lambda K,\cdot)
	=
	\lambda^{(n+q-1)s-1}F_{s,q}(K,\cdot).
	\]
	Consequently, for $K\in\Kno$ and $p\in\bR$,
	\[
	F_{p,s,q}(\lambda K,\cdot)
	=
	\lambda^{(n+q-1)s-p}F_{p,s,q}(K,\cdot).
	\]
\end{proposition}

\begin{proof}
	Let $\omega\subset\Sn$ be a Borel set. Since
	\[
	\nu_{\lambda K}(\lambda z)=\nu_K(z)
	\]
	for $\cHnmone$-almost every $z\in\partial K$, and
	\[
	d\cHnmone(\lambda z)=\lambda^{n-1}d\cHnmone(z),
	\]
	Definition~\ref{Def-sq-chord-measure} gives
	\begin{align*}
		F_{s,q}(\lambda K,\omega)
		&=
		\frac{q}{n\omega_n}
		\int_{\nu_K^{-1}(\omega)}
		\int_{\Sn}
		I_q(\lambda K,u)^{s-1}
		X_{\lambda K}(\lambda z,u)^{q-1}
		\,du\,
		\lambda^{n-1}d\cHnmone(z).
	\end{align*}
	By \eqref{eq-homogeneity-of-Iq} and
	\[
	X_{\lambda K}(\lambda z,u)=\lambda X_K(z,u),
	\]
	we have
	\begin{align*}
		F_{s,q}(\lambda K,\omega)
		=
		\lambda^{(n+q-1)s-1}F_{s,q}(K,\omega).
	\end{align*}
	Since $\omega$ is arbitrary,
	\[
	F_{s,q}(\lambda K,\cdot)
	=
	\lambda^{(n+q-1)s-1}F_{s,q}(K,\cdot).
	\]
	
	Finally, since $h_{\lambda K}=\lambda h_K$, we obtain
	\[
	F_{p,s,q}(\lambda K,\cdot)
	=
	\lambda^{(n+q-1)s-p}F_{p,s,q}(K,\cdot).
	\] The assertions follow.
\end{proof}

\begin{proposition}\label{prop-evenness-of-sq-chord-measure}
	If $K\in\Kne$, then $F_{s,q}(K,\cdot)$ is even for every $q>0$ and $s\ne0$.
\end{proposition}

\begin{proof}
	Let $\omega\subset\Sn$ be Borel. Since $K=-K$,
	\[
	\nu_K(-y)=-\nu_K(y)
	\quad\text{and}\quad
	X_K(-y,u)=X_K(y,u)
	\]
	for $\cHnmone$-almost every $y\in\partial K$. Hence
	\[
	\nu_K^{-1}(-\omega)=-\nu_K^{-1}(\omega).
	\]
	Using Definition~\ref{Def-sq-chord-measure} and the change of variables
	$y=-z$, we obtain
	\begin{align*}
		F_{s,q}(K,-\omega)
		&=
		\frac{q}{n\omega_n}
		\int_{\nu_K^{-1}(\omega)}
		\int_{\Sn}
		I_q(K,u)^{s-1}X_K(-z,u)^{q-1}
		\,du\,d\cHnmone(z)\\
		&=
		F_{s,q}(K,\omega).
	\end{align*}
	Therefore $F_{s,q}(K,\cdot)$ is even.
\end{proof}

Now we formulate the $L_p$ Minkowski problem for $(s,q)$-chord measures. Given $q>0$, $s\ne0$, $p\in\bR$, and a finite Borel measure $\mu$ on $\Sn$, determine when there exists a convex body $K$ such that
\[
F_{p,s,q}(K,\cdot)=\mu,
\]
with $K\in\Kno$.

Suppose that $d\mu(v)=f(v)\,dv$. We need to find a positive function $h\in C^2(\Sn)$ such that for every $v\in\mathbb S^{n-1}$,
\[
\nabla_S^2h(v)+h(v)I>0
\qquad\text{on }v^\perp,
\]
where $\nabla_S^2h$ is the spherical Hessian and $I$ is the identity matrix on $v^\perp$. Let $K_h=[h]$ and
\[
x_h(v)=\nabla_Sh(v)+h(v)v\in\partial K_h.
\]
Then the measure equation $F_{p,s,q}(K_h,\cdot)=\mu$ becomes the nonlocal Monge--Amp\`ere equation
\begin{equation}\label{eqsec3-MA-equation}
	\frac{q}{n\omega_n}
	h(v)^{1-p}
	\left[
	\int_{\Sn}
	I_q(K_h,u)^{s-1}
	X_{K_h}\bigl(x_h(v),u\bigr)^{q-1}\,du
	\right]
	\det\bigl(\nabla_S^2h(v)+h(v)I\bigr)
	=f(v).
\end{equation}

\section{\texorpdfstring{Minkowski problems for $(s,q)$-chord measures}{Minkowski problems for (s,q)-chord measures}}\label{sec:minkowski}

For a nonzero finite Borel measure $\mu$ on $\Sn$ and $f\in C^+(\Sn)$, write
\[
\langle f,\mu\rangle=\int_{\Sn}f\,d\mu.
\]
Define the $0$-homogeneous functional $\Phi_{s,q,\mu}:\Kno\to(0,\infty)$ by
\[
\Phi_{s,q,\mu}(K)
=
\frac{\Psi_{s,q}(K)^{1/((n+q-1)s)}}{\langle h_K,\mu\rangle}.
\]

\begin{lemma}\label{lem-characterization-of-maximizer}
Let $q>0$ and $s\in\bR\setminus\{0\}$ satisfy $(n+q-1)s\ne1$, and let $\mu$ be a nonzero finite Borel measure on $\Sn$. If $\Phi_{s,q,\mu}$ has a maximizer in $\Kno$, then there
exists $K\in\Kno$ such that
\[
F_{s,q}(K,\cdot)=\mu.
\]
\end{lemma}

\begin{proof}
Let $K_0\in\Kno$ be a maximizer of $\Phi_{s,q,\mu}$. For $g\in C(\Sn)$, set $h_t=h_{K_0}+tg$ and $K_t=[h_t]$ for sufficiently small $|t|$, so that $h_t>0$. Define
\[
\varphi(t)
=
\frac{\Psi_{s,q}(K_t)^{1/((n+q-1)s)}}{\langle h_t,\mu\rangle}.
\]
Since $h_{K_t}\le h_t$ and $\mu$ is nonnegative,
\[
\varphi(t)
\le
\Phi_{s,q,\mu}(K_t)
\le
\Phi_{s,q,\mu}(K_0)
=
\varphi(0).
\]
Thus $\varphi'(0)=0$. Theorem~\ref{thm-variation-formula-for-sq-chord-int} gives
\begin{align*}
0
&=
\frac{\Psi_{s,q}(K_0)^{(1-(n+q-1)s)/((n+q-1)s)}}{(n+q-1)\langle h_{K_0},\mu\rangle}
\int_{\Sn}g(v)\,dF_{s,q}(K_0,v)\\
&\quad-
\frac{\Psi_{s,q}(K_0)^{1/((n+q-1)s)}}{\langle h_{K_0},\mu\rangle^2}
\int_{\Sn}g(v)\,d\mu(v).
\end{align*}
Since this holds for every $g\in C(\Sn)$,
\[
\frac{\langle h_{K_0},\mu\rangle}{(n+q-1)\Psi_{s,q}(K_0)}
F_{s,q}(K_0,\cdot)
=
\mu.
\]
Let $\lambda>0$ be such that
\[
\lambda^{(n+q-1)s-1}
=
\frac{\langle h_{K_0},\mu\rangle}{(n+q-1)\Psi_{s,q}(K_0)}.
\] Let $K=\lambda K_0$. Proposition~\ref{prop-homogeneity-of-sq-chord-measure} then shows that $K$ satisfies $F_{s,q}(K,\cdot)=\mu$.
\end{proof}

\begin{proof}[Proof of Theorem~\ref{thm-mp-for-sq-chord-measure}]
Assume that $\mu$ is not concentrated on any closed hemisphere and satisfies \eqref{eq-centroid-of-mu}. Consider
\[
\sup\{\Phi_{s,q,\mu}(K):K\in\Kno\}.
\]
Because $\int_{\Sn}v\,d\mu(v)=o$ and $\Psi_{s,q}$ is translation invariant, $\Phi_{s,q,\mu}$ is also translation invariant. We may therefore choose a maximizing sequence $(K_i)$ of \(\Phi_{s,q,\mu}\) such that $c(K_i)=o$ for every $i$. By $0$-homogeneity of \(\Phi_{s,q,\mu}\), we may also normalize \((K_i)\) so that
\[
\Psi_{s,q}(K_i)=1.
\]
Under this normalization, $(K_i)$ is a minimizing sequence for
\begin{equation*}
\inf\left\{
\langle h_K,\mu\rangle:
K\in\Kno,\ \Psi_{s,q}(K)=1
\right\}.
\end{equation*}
Choose $a>0$ such that $\Psi_{s,q}(aB^n)=1$. After discarding finitely many terms of \((K_i)\), we can assume that
\[
	\langle h_{K_i},\mu\rangle
	\le
	\langle h_{aB^n},\mu\rangle+1
	=a|\mu|+1.
\]

Let
\[
R_i=\max_{u\in\Sn}\rho_{K_i}(u)=\rho_{K_i}(u_i)
\]
for some $u_i\in\Sn$. Since $[o,R_iu_i]\subset K_i$,
\[
R_i(v\cdot u_i)_+
=h_{[o,R_iu_i]}(v)
\le h_{K_i}(v),
\]
where $(t)_+=\max\{t,0\}$. Hence
\[
R_i\int_{\Sn}(u_i\cdot v)_+\,d\mu(v)
\le a|\mu|+1.
\]
The continuous function
\[
u\longmapsto\int_{\Sn}(u\cdot v)_+\,d\mu(v)
\]
is strictly positive on $\Sn$ because $\mu$ is not concentrated on a closed hemisphere. Therefore there is a constant $c>0$ such that
\[
\int_{\Sn}(u\cdot v)_+\,d\mu(v)\ge c
\qquad\text{for every }u\in\Sn.
\]
Thus $R_i\le (a|\mu|+1)/c$, so $(K_i)$ is uniformly bounded. By the Blaschke selection theorem, a subsequence, which is denoted by $(K_i)$ again, converges in the Hausdorff metric to a compact convex set $K_0$.

We claim that $K_0$ has nonempty interior. Otherwise, there is $e\in\Sn$ such that $w_{K_0}(e)=0$, and hence $w_{K_i}(e)\to0$. For $u\cdot e\ne0$, every chord of $K_i$ parallel to $u$ satisfies
\[
X_{K_i}(x,u)\le\frac{w_{K_i}(e)}{|u\cdot e|}.
\]
Since $K_i\subset RB^n$ for some $R>0$ independent of $i$,
\begin{align*}
I_q(K_i,u)
&\le
\vol_{n-1}(K_i|u^\perp)
\left(\frac{w_{K_i}(e)}{|u\cdot e|}\right)^q\\
&\le
\omega_{n-1}R^{n-1}
\left(\frac{w_{K_i}(e)}{|u\cdot e|}\right)^q.
\end{align*}
Thus $I_q(K_i,u)\to0$ for almost every $u\in\Sn$. For $s>0$, since
\[
I_q(K_i,u)
\le
\omega_{n-1}R^{n-1}(2R)^q,
\]
 dominated convergence gives $\Psi_{s,q}(K_i)\to0$, a contradiction. If $s<0$, then $I_q(K_i,u)^s\to+\infty$ for almost every $u$, and Fatou's lemma gives $\liminf_i\Psi_{s,q}(K_i)=+\infty$, again a contradiction. Hence $K_0\in\Kn$.

Hausdorff continuity of $\Psi_{s,q}$ gives $\Psi_{s,q}(K_0)=1$. Since centroids are continuous under the Hausdorff metric, $c(K_0)=o$; in particular, $K_0\in\Kno$. Because $\Phi_{s,q,\mu}$ is continuous under the Hausdorff metric
and $(K_i)$ is a maximizing sequence for $\Phi_{s,q,\mu}$, we have \[\Phi_{s,q,\mu}(K_0)=\sup\{\Phi_{s,q,\mu}(K):K\in\Kno\}.\] Lemma~\ref{lem-characterization-of-maximizer} therefore yields a dilation $K$ of $K_0$ satisfying
\[
F_{s,q}(K,\cdot)=\mu.
\]

Conversely, suppose that $F_{s,q}(K,\cdot)=\mu$ for some $K\in\Kn$, and define
\[
B_K(z)
=
\frac{q}{n\omega_n}
\int_{\Sn}I_q(K,u)^{s-1}X_K(z,u)^{q-1}\,du,
\qquad z\in\partial K.
\]
Then
\[
F_{s,q}(K,\omega)
=
\int_{\nu_K^{-1}(\omega)}B_K(z)\,d\cHnmone(z)
\]
for every Borel set $\omega\subset\Sn$. Lemma~\ref{lemsec2-finite-of-Fsq} gives $B_K(z)<\infty$ for $\cHnmone$-almost every $z$. To prove positivity of $B_K(z)$, fix $z\in\partial K$, choose $a\in\operatorname{int}K$ and $r>0$ with $B(a,r)=\{x\in\mathbb R^n:|x-a|<r\}\subset K$, and let
\[
u_0=\frac{a-z}{|a-z|}.
\]
For every $u$ in a sufficiently small spherical neighborhood $U$ of $u_0$, the line $z+\mathbb Ru$ meets $B(a,r)$ in a nondegenerate segment. Hence $X_K(z,u)>0$ on $U$, and therefore $B_K(z)>0$ whenever it is finite. Consequently, if $F_{s,q}(K,\omega)=0$
then
$S(K,\omega)=0$.
Thus $S(K,\cdot)$ is absolutely continuous with respect to $F_{s,q}(K,\cdot)$. Since a surface area measure is not concentrated on any closed hemisphere, neither is $F_{s,q}(K,\cdot)$. Finally, Theorem~\ref{thm-centroid-of-sq-chord-measure} gives
\[
\int_{\Sn}v\,dF_{s,q}(K,v)=o.
\]
This proves the necessity and completes the proof.
\end{proof}

We next establish the inequality used for uniqueness. The classical Brunn--Minkowski inequality shows that, for $K,L\in\Kn$ and $0<\lambda<1$,
\begin{equation}\label{eq-first-ClassicalBm}
V((1-\lambda)K+\lambda L)
\ge
V(K)^{1-\lambda}V(L)^\lambda,
\end{equation}
with equality if and only if $K$ and $L$ are translates; see \cite{schneiderbook}.

Let $u\in\mathbb S^{n-1}$. For $y\in K|u^\perp$, define the upper and lower graph functions in direction $u$ by
\[
\overline\ell_u(K,y)=\max\{t\in\bR:y+tu\in K\},
\qquad
\underline\ell_u(K,y)=-\min\{t\in\bR:y+tu\in K\}.
\]
Then
\[
K=\{y+tu:y\in K|u^\perp,\ -\underline\ell_u(K,y)\le t\le\overline\ell_u(K,y)\},
\]
and both graph functions are concave on $K|u^\perp$.  For  $y\in K|u^\perp$, we have
\begin{equation}\label{eq:Xgraph}
X_K(y,u)=\overline\ell_u(K,y)+\underline\ell_u(K,y).
\end{equation}

We will use the Pr\'ekopa--Leindler inequality together with its equality characterization.

\begin{theorem}[Pr\'ekopa--Leindler--Dubuc \cite{Dubuc1977}]\label{thm:PL}
Let $0<\lambda<1$, and let $f,g,h:\bR^m\to[0,\infty)$ be integrable functions satisfying
\[
h((1-\lambda)x+\lambda y)
\ge
f(x)^{1-\lambda}g(y)^\lambda
\]
for all $x,y\in\bR^m$. Then
\[
\int_{\bR^m}h
\ge
\left(\int_{\bR^m}f\right)^{1-\lambda}
\left(\int_{\bR^m}g\right)^\lambda.
\]
Assume in addition that
$\int_{\bR^m}f>0$ and $\int_{\bR^m}g>0$.
Then equality holds if and only if there exist $a>0$,
$z\in\bR^m$, and a log-concave function $H$ such that,
almost everywhere,
\[
h=H,\qquad
f(x)=a^{-\lambda}H(x-\lambda z),\qquad
g(x)=a^{1-\lambda}H(x+(1-\lambda)z).
\]
\end{theorem}

\begin{lemma}\label{lem:XrayBM}
Let $K,L\in\Kn$, $0<\lambda<1$, $u\in\Sn$, and $K_\lambda=(1-\lambda)K+\lambda L$. Then, for all $x,y\in u^\perp$,
\[
X_{K_\lambda}((1-\lambda)x+\lambda y,u)
\ge
X_K(x,u)^{1-\lambda}X_L(y,u)^\lambda.
\]
\end{lemma}

\begin{proof}
	If $x\notin K|u^\perp$ or $y\notin L|u^\perp$, then the right-hand side is
	zero, and the assertion is immediate. Hence assume
	$x\in K|u^\perp$ and $y\in L|u^\perp$.
	
	Write
	\[
	X_K(x,u)=\bar\ell_u(K,x)+\underline\ell_u(K,x),
	\qquad
	X_L(y,u)=\bar\ell_u(L,y)+\underline\ell_u(L,y).
	\]
	 Since
	\[
	x+\bar\ell_u(K,x)u\in K,
	\qquad
	y+\bar\ell_u(L,y)u\in L,
	\]
	we have
	\[
	(1-\lambda)x+\lambda y+
	\bigl((1-\lambda)\bar\ell_u(K,x)
	+\lambda\bar\ell_u(L,y)\bigr)u
	\in K_\lambda.
	\]
	Thus
	\[
	\bar\ell_u(K_\lambda,(1-\lambda)x+\lambda y)
	\ge
	(1-\lambda)\bar\ell_u(K,x)
	+\lambda\bar\ell_u(L,y).
	\]
	Similarly,
	\[
	\underline\ell_u(K_\lambda,(1-\lambda)x+\lambda y)
	\ge
	(1-\lambda)\underline\ell_u(K,x)
	+\lambda\underline\ell_u(L,y).
	\]
	Adding these inequalities gives
	\[
	X_{K_\lambda}((1-\lambda)x+\lambda y,u)
	\ge
	(1-\lambda)X_K(x,u)+\lambda X_L(y,u).
	\]
	Finally, by the weighted arithmetic--geometric mean inequality,
	\[
	(1-\lambda)X_K(x,u)+\lambda X_L(y,u)
	\ge
	X_K(x,u)^{1-\lambda}X_L(y,u)^\lambda.
	\]
	This proves the lemma.
\end{proof}

\begin{lemma}\label{lem:Xlogconcave}
For every $K\in\Kn$ and $u\in\Sn$, the function $x\mapsto X_K(x,u)$ is log-concave on $u^\perp$.
\end{lemma}

\begin{proof}
	Let $x,y\in u^\perp$ and $0<\lambda<1$. If either
	$x\notin K|u^\perp$ or $y\notin K|u^\perp$, then
	\[
	X_K(x,u)^{1-\lambda}X_K(y,u)^\lambda=0,
	\]
	so the log-concavity inequality is immediate.
	
	Now assume $x,y\in K|u^\perp$. By \eqref{eq:Xgraph}, we have
	\begin{align*}
		X_K((1-\lambda)x+\lambda y,u)
		&\ge
		(1-\lambda)X_K(x,u)+\lambda X_K(y,u)\\
		&\ge
		X_K(x,u)^{1-\lambda}X_K(y,u)^\lambda,
	\end{align*}
	where the second inequality is the weighted arithmetic--geometric mean inequality.
	
	Therefore $x\mapsto X_K(x,u)$ is log-concave on $u^\perp$.
\end{proof}

\begin{theorem}\label{thm:firstBM}
Let $K,L\in\Kn$, $s<0$, $q>0$, and $0<\lambda<1$. Then
\begin{equation}\label{eq:firstBM}
\Psi_{s,q}((1-\lambda)K+\lambda L)
\le
\Psi_{s,q}(K)^{1-\lambda}\Psi_{s,q}(L)^\lambda.
\end{equation}
Equality holds if and only if $K$ and $L$ are translates.
\end{theorem}

\begin{proof}
Set $K_\lambda=(1-\lambda)K+\lambda L$ and fix $u\in\Sn$. By Lemma~\ref{lem:XrayBM} and Theorem~\ref{thm:PL}, we have
\begin{equation}\label{eq:JBM}
I_q(K_\lambda,u)
\ge
I_q(K,u)^{1-\lambda}I_q(L,u)^\lambda.
\end{equation}
Since $s<0$,
\[
I_q(K_\lambda,u)^s
\le
I_q(K,u)^{(1-\lambda)s}I_q(L,u)^{\lambda s}.
\]Integrating the above inequality, we have
\begin{align*}
	\Psi_{s,q}(K_\lambda)
	&=\frac{1}{n\omega_n}\int_{\Sn}I_q(K_\lambda,u)^s\,du\\
	&\le
	\frac{1}{n\omega_n}\int_{\Sn}
	I_q(K,u)^{(1-\lambda)s}
	I_q(L,u)^{\lambda s}\,du\\
	&\le
	\left(
	\frac{1}{n\omega_n}\int_{\Sn}I_q(K,u)^s\,du
	\right)^{1-\lambda}
	\left(
	\frac{1}{n\omega_n}\int_{\Sn}I_q(L,u)^s\,du
	\right)^\lambda\\
	&=
	\Psi_{s,q}(K)^{1-\lambda}
	\Psi_{s,q}(L)^\lambda.
\end{align*}
Thus we prove the desired inequality.

If $K$ and $L$ are translates, equality follows from translation invariance of \(\Psi_{s,q}\). Conversely, suppose that equality holds. This implies equality holds in \eqref{eq:JBM} for almost every $u$. Fix such a direction. Apply the equality statement in Theorem~\ref{thm:PL} to
\[
f(x)=X_K(x,u)^q,\quad
g(x)=X_L(x,u)^q,\quad
h(x)=X_{K_\lambda}(x,u)^q.
\]
Lemma~\ref{lem:Xlogconcave} shows that $h$ is log-concave. Hence there are $a>0$ and $z\in u^\perp$ such that, almost everywhere on $u^\perp$,
\[
X_K(x,u)
=
a^{-\lambda/q}X_{K_\lambda}(x-\lambda z,u),
\]
and
\[
X_L(x,u)
=
a^{(1-\lambda)/q}X_{K_\lambda}(x+(1-\lambda)z,u).
\]
Integrating over $u^\perp$ gives
\[
V(K)=a^{-\lambda/q}V(K_\lambda),
\qquad
V(L)=a^{(1-\lambda)/q}V(K_\lambda).
\]
Therefore
\[
V(K_\lambda)=V(K)^{1-\lambda}V(L)^\lambda.
\]
Equality in \eqref{eq-first-ClassicalBm} implies that $K$ and $L$ are translates.
\end{proof}

\begin{theorem}\label{thm:secondBM}
Let $K,L\in\Kn$, $s<0$, and $q>0$. Then
\begin{equation}\label{eq:secondBM}
\Psi_{s,q}(K+L)^{1/((n+q-1)s)}
\ge
\Psi_{s,q}(K)^{1/((n+q-1)s)}
+
\Psi_{s,q}(L)^{1/((n+q-1)s)}.
\end{equation}
Equality holds if and only if $K$ and $L$ are homothetic.
\end{theorem}

\begin{proof}
Let $\alpha=(n+q-1)s<0$ and set
\[
a=\Psi_{s,q}(K)^{1/\alpha},
\qquad
b=\Psi_{s,q}(L)^{1/\alpha},
\qquad
\widetilde K=K/a,
\qquad
\widetilde L=L/b.
\]
By Lemma~\ref{lem-homogeneity-of-sq-chord-int}, $\Psi_{s,q}(\widetilde K)=\Psi_{s,q}(\widetilde L)=1$. With $\lambda=b/(a+b)$,
\[
(1-\lambda)\widetilde K+\lambda\widetilde L
=\frac{K+L}{a+b}.
\]
Theorem~\ref{thm:firstBM} gives
\[
\Psi_{s,q}\left(\frac{K+L}{a+b}\right)\le1.
\]
Since $1/\alpha<0$, homogeneity yields \eqref{eq:secondBM}.

Equality in \eqref{eq:secondBM} is equivalent to equality in Theorem~\ref{thm:firstBM} for $\widetilde K$ and $\widetilde L$, and hence $\widetilde{K},\widetilde{L}$ are translates. This implies that $K$ and $L$ are homothetic. The converse follows directly from translation invariance and homogeneity.
\end{proof}
\begin{theorem}\label{thm-sq-chord-measure-uniquely-determine-body-when-s-le-0}
Let $q>0$, $s<0$, and $K,L\in\Kn$. If
\[
F_{s,q}(K,\cdot)=F_{s,q}(L,\cdot),
\]
then $K$ and $L$ are translates.
\end{theorem}

\begin{proof}
Let $K_t=K+tL$ for $t\ge0$ and put $\alpha=(n+q-1)s<0$. By Theorems~\ref{thm-variation-formula-for-sq-chord-int} and \ref{thm-integral-representation-of-sq-chord},
\begin{align}
\left.\frac{d}{dt}\right|_{t=0^+}\Psi_{s,q}(K_t)^{1/\alpha}
&=
\frac{1}{n+q-1}\Psi_{s,q}(K)^{(1-\alpha)/\alpha}
\int_{\Sn}h_L(v)\,dF_{s,q}(K,v)\notag\\
&=
\frac{1}{n+q-1}\Psi_{s,q}(K)^{(1-\alpha)/\alpha}
\int_{\Sn}h_L(v)\,dF_{s,q}(L,v)\notag\\
&=
\Psi_{s,q}(K)^{(1-\alpha)/\alpha}\Psi_{s,q}(L).
\label{eq-mixed-sq-chord-for-K}
\end{align}
On the other hand, Theorem~\ref{thm:secondBM} gives
\begin{equation}\label{eq-limt-BM-ineq}
\left.\frac{d}{dt}\right|_{t=0^+}\Psi_{s,q}(K_t)^{1/\alpha}
\ge
\Psi_{s,q}(L)^{1/\alpha}.
\end{equation}
Since $(1-\alpha)/\alpha<0$, \eqref{eq-mixed-sq-chord-for-K} and \eqref{eq-limt-BM-ineq} imply $\Psi_{s,q}(K)\le\Psi_{s,q}(L)$. Interchanging the roles of $K$ and $L$ gives
\begin{equation}\label{eq-Psi-sq-K-Psi-sq-L}
\Psi_{s,q}(K)=\Psi_{s,q}(L).
\end{equation}

Define
\[
f(t)=\Psi_{s,q}(K+tL)^{1/\alpha}
-\Psi_{s,q}(K)^{1/\alpha}
-t\Psi_{s,q}(L)^{1/\alpha}.
\]
Theorem~\ref{thm:secondBM} gives $f(t)\ge0$, and $f$ is concave on $[0,\infty)$. \eqref{eq-mixed-sq-chord-for-K} and \eqref{eq-Psi-sq-K-Psi-sq-L} give $f(0)=f'(0^+)=0$. Concavity implies $f\le0$, and hence $f\equiv0$. Equality in Theorem~\ref{thm:secondBM} shows that $K=\lambda L+x_0$ for $\lambda>0$ and $x_0\in\mathbb R^n$. \eqref{eq-Psi-sq-K-Psi-sq-L} shows that $\lambda=1$, so $K$ and $L$ are translates.
\end{proof}
Using Theorem \ref{thm-mp-for-sq-chord-measure} and Theorem \ref{thm-sq-chord-measure-uniquely-determine-body-when-s-le-0}, we have the following corollary.
\begin{corollary}\label{coro-unique-of-mp}
Let $q>0$, $s<0$, and let $\mu$ be a nonzero finite Borel measure on $\Sn$. Then there exists a unique convex body, up to translation, satisfying
\[
F_{s,q}(K,\cdot)=\mu
\]
if and only if $\mu$ is not concentrated on any closed hemisphere and
\[
\int_{\Sn}v\,d\mu(v)=o.
\]
\end{corollary}

\section{\texorpdfstring{The logarithmic Minkowski problem for $(s,q)$-chord measures}{The logarithmic Minkowski problem for (s,q)-chord measures}}\label{sec:log-minkowski}

In this section we study the logarithmic Minkowski problem in the origin-symmetric case. Let $\mu$ be a nonzero finite Borel measure on $\Sn$. For $K\in\Kno$, define the entropy
\[
E_\mu(K)
=-\frac{1}{|\mu|}\int_{\Sn}\log h_K(v)\,d\mu(v).
\]

Let $q\ge1$ and $s>0$. 
We say that $\mu$ satisfies the subspace mass inequality if
\begin{equation}\label{eqsec4-subspace-mass-inequality}
\frac{\mu(\xi\cap\Sn)}{|\mu|}
<
\frac{ks+\min\{k,(q-1)s\}}{(n+q-1)s}
\end{equation}
for every $k\in\{1,\ldots,n-1\}$ and every $\xi\in G_{n,k}$.

Let $e_1,\ldots,e_n$ be an orthonormal basis of $\Rn$. For $0<\delta<1/\sqrt n$, define
\begin{equation}\label{eqsec4-partition-of-sphere}
\Omega_{k,\delta}
=
\left\{v\in\Sn:
|v\cdot e_k|\ge\delta,
\ |v\cdot e_j|<\delta\text{ for all }j>k
\right\},
\qquad k=1,\ldots,n.
\end{equation}
These sets form a partition of $\Sn$. Put $\xi_0=\{o\}$ and
\[
\xi_k=\operatorname{span}\{e_1,\ldots,e_k\},
\qquad k=1,\ldots,n.
\]
It was proved in \cite{Boroczky2013JAMS} that
\begin{equation}\label{eqsec4-lim-of-partition-measure}
\lim_{\delta\to0^+}\mu(\Omega_{k,\delta})
=
\mu\bigl((\xi_k\setminus\xi_{k-1})\cap\Sn\bigr),
\end{equation}
and therefore
\begin{equation}\label{eqsec4-lim-sum-of-partition-measure}
\lim_{\delta\to0^+}
\sum_{j=1}^k\mu(\Omega_{j,\delta})
=
\mu(\xi_k\cap\Sn).
\end{equation}

We also use the following lemma.

\begin{lemma}[Lemma~4.1, \cite{Boroczky2019Adv}]\label{lemsec4-sum-by-part-ineq}
Let $\lambda_1,\ldots,\lambda_m\in[0,1]$ satisfy $\sum_{j=1}^m\lambda_j=1$, and let $a_1\le\cdots\le a_m$. Suppose that $0=\sigma_0\le\sigma_1\le\cdots\le\sigma_m=1$ and
\[
\sum_{j=1}^k\lambda_j\le\sigma_k,
\qquad k=1,\ldots,m.
\]
Then
\[
\sum_{k=1}^m\lambda_k a_k
\ge
\sum_{k=1}^m(\sigma_k-\sigma_{k-1})a_k.
\]
\end{lemma}
For $\beta\in\mathbb R$, let $\lceil\beta\rceil$ denote the smallest integer greater than or equal to \(\beta\). We first establish the following entropy estimate. 
\begin{lemma}\label{lemsec4-entropy-estimate}
Let $q\ge1$ and $s>0$, set $\beta=(q-1)s$, assume $\beta\in[0,n]$, and put $m=\max\{1,\lceil\beta\rceil\}$. Let $\mu$ be a nonzero finite Borel measure satisfying \eqref{eqsec4-subspace-mass-inequality}. Suppose that the ordered orthonormal bases $(e_{1l},\ldots,e_{nl})$ converge to $(e_1,\ldots,e_n)$. Let $\epsilon_0>0$, and let
\[
E_l
=
\left\{x\in\Rn:
\sum_{k=1}^n\frac{(x\cdot e_{kl})^2}{a_{k,l}^2}\le1
\right\},
\]
where
\[
0<a_{1,l}\le\cdots\le a_{n,l},
\qquad a_{n,l}\ge\epsilon_0.
\]
Then there exist $l_0\in\mathbb N$, $t_0>0$, and $c_0\in\bR$ such that, for every $l\ge l_0$,
\begin{align}
E_\mu(E_l)
\le{}&
-\frac{s+1}{(n+q-1)s}\log\left(\prod_{k=1}^{m-1}a_{k,l}\right)
-\frac{\beta+s+1-m}{(n+q-1)s}\log a_{m,l}\notag\\
&-\frac{s}{(n+q-1)s}\log\left(\prod_{k=m+1}^{n}a_{k,l}\right)
+t_0\log a_{1,l}+c_0.
\label{eqsec4-entropy-estimate-for-ellipsoid}
\end{align}
Empty products are interpreted as $1$.
\end{lemma}

\begin{proof}
	For the partition \eqref{eqsec4-partition-of-sphere}, let
	\[
	\lambda_{k,\delta}
	=
	\frac{\mu(\Omega_{k,\delta})}{|\mu|},
	\qquad k=1,\ldots,n.
	\]
	By \eqref{eqsec4-lim-sum-of-partition-measure} and
	\eqref{eqsec4-subspace-mass-inequality}, there exist
	$t_0>0$ and $\delta_0\in(0,1/\sqrt n)$ such that
	\[
	\sum_{j=1}^k\lambda_{j,\delta_0}
	<
	\frac{ks+\min\{k,\beta\}}{(n+q-1)s}-t_0,
	\qquad k=1,\ldots,n-1.
	\]
	 Define
	\[
	\sigma_0=0,\qquad
	\sigma_k=\frac{ks+\min\{k,\beta\}}{(n+q-1)s}-t_0
	\quad (1\le k\le n-1),
	\qquad
	\sigma_n=1.
	\]
	Let $t_0$ be sufficiently small. Then $0=\sigma_0\le\sigma_1\le\cdots\le\sigma_n=1$.
	
	 Since $m=\max\{1,\lceil\beta\rceil\}$, a direct calculation gives
	\[
	\sigma_1-\sigma_0
	=
	\frac{\min\{s+1,s+\beta\}}{(n+q-1)s}-t_0,
	\]
	and, for $2\le k\le n-1$,
	\[
	\sigma_k-\sigma_{k-1}
	=
	\begin{cases}
		\dfrac{s+1}{(n+q-1)s}, & k<m,\\[4pt]
		\dfrac{\beta+s+1-m}{(n+q-1)s}, & k=m,\\[4pt]
		\dfrac{s}{(n+q-1)s}, & k>m.
	\end{cases}
	\]
	Moreover,
	\[
	\sigma_n-\sigma_{n-1}
	=
	\frac{s+\beta-\min\{n-1,\beta\}}{(n+q-1)s}+t_0.
	\]
These formulas also hold when $m=1$ or $m=n$.
	
	Since $(e_{1l},\ldots,e_{nl})\to(e_1,\ldots,e_n)$, there exists
	$l_0\in\mathbb N$ such that, for every $l\ge l_0$,
	\[
	|e_{kl}-e_k|<\frac{\delta_0}{2},
	\qquad k=1,\ldots,n.
	\]
	For $v\in\Omega_{k,\delta_0}$, since $a_{k,l}e_{kl}\in E_l$, we obtain
	\begin{align*}
		h_{E_l}(v)
		&\ge |v\cdot a_{k,l}e_{kl}|\\
		&\ge a_{k,l}
		\bigl(|v\cdot e_k|-|v\cdot(e_{kl}-e_k)|\bigr)\\
		&\ge \frac{\delta_0}{2}a_{k,l}.
	\end{align*}
	Hence
	\begin{align}
		E_\mu(E_l)
		&=
		-\frac1{|\mu|}
		\int_{\Sn}\log h_{E_l}(v)\,d\mu(v)\notag\\
		&\le
		-\sum_{k=1}^n
		\lambda_{k,\delta_0}
		\log\left(\frac{\delta_0}{2}a_{k,l}\right)\notag\\
		&=
		-\log\frac{\delta_0}{2}
		-\sum_{k=1}^n\lambda_{k,\delta_0}\log a_{k,l}\notag.
	\end{align}
	Since
	\[
	\log a_{1,l}\le\cdots\le\log a_{n,l},
	\]
	Lemma~\ref{lemsec4-sum-by-part-ineq} yields
	\[
	\sum_{k=1}^n\lambda_{k,\delta_0}\log a_{k,l}
	\ge
	\sum_{k=1}^n
	(\sigma_k-\sigma_{k-1})\log a_{k,l}.
	\]
	Therefore
	\begin{align*}
		E_\mu(E_l)
		\le{}&
		-\log\frac{\delta_0}{2}
		-\sum_{k=1}^n
		(\sigma_k-\sigma_{k-1})\log a_{k,l}.
	\end{align*}
	
	Substituting the above formulas for $\sigma_k-\sigma_{k-1}$, we get
	\begin{align*}
		-\sum_{k=1}^n
		(\sigma_k-\sigma_{k-1})\log a_{k,l}
		={}&
		-\frac{s+1}{(n+q-1)s}
		\sum_{k=1}^{m-1}\log a_{k,l}\\
		&-\frac{\beta+s+1-m}{(n+q-1)s}\log a_{m,l}\\
		&-\frac{s}{(n+q-1)s}
		\sum_{k=m+1}^n\log a_{k,l}\\
		&+t_0\log a_{1,l}-t_0\log a_{n,l}.
	\end{align*}
	Consequently,
	\begin{align*}
		E_\mu(E_l)
		\le{}&
		-\frac{s+1}{(n+q-1)s}
		\log\left(\prod_{k=1}^{m-1}a_{k,l}\right)
		-\frac{\beta+s+1-m}{(n+q-1)s}\log a_{m,l}\\
		&-\frac{s}{(n+q-1)s}
		\log\left(\prod_{k=m+1}^n a_{k,l}\right)
		+t_0\log a_{1,l}
		-\!t_0\log a_{n,l}
		-\log\frac{\delta_0}{2}.
	\end{align*}
	Finally, since $a_{n,l}\ge\epsilon_0$,
	\[
	-t_0\log a_{n,l}
	\le
	-t_0\log\epsilon_0.
	\]
	Let \[c_0=-t_0\log\epsilon_0-\log\frac{\delta_0}2.\] We have
	\begin{align*}
		E_\mu(E_l)
		\le{}&
		-\frac{s+1}{(n+q-1)s}
		\log\left(\prod_{k=1}^{m-1}a_{k,l}\right)-\frac{\beta+s+1-m}{(n+q-1)s}\log a_{m,l}\\
		&-\frac{s}{(n+q-1)s}
		\log\left(\prod_{k=m+1}^n a_{k,l}\right)
		+t_0\log a_{1,l}+c_0.
	\end{align*}
	This is \eqref{eqsec4-entropy-estimate-for-ellipsoid}.
\end{proof}

We next give an $(s,q)$-chord estimate for ellipsoids.

\begin{lemma}\label{lemsec4-ellipsoid-estimate}
Let $q>1$ and $s>0$, and set $\beta=(q-1)s$.
Assume that $\beta\in(0,n)\setminus\mathbb Z$, and put
$m=\lceil\beta\rceil$. Let $e_1,\ldots,e_n$ be an orthonormal basis and
\begin{equation}\label{eqsec4-ellipsoid}
E
=
\left\{x\in\Rn:
\sum_{k=1}^n\frac{(x\cdot e_k)^2}{a_k^2}\le1
\right\},
\qquad
0<a_1\le\cdots\le a_n\le1.
\end{equation}
Then
\[
\Psi_{s,q}(E)
\le
c_{s,q,n}
\left(\prod_{k=1}^{m-1}a_k\right)^{s+1}
 a_m^{\beta+s+1-m}
\left(\prod_{k=m+1}^{n}a_k\right)^s,
\]where $c_{s,q,n}$ is a constant depending only on $s,q,n$.
\end{lemma}

\begin{proof}
Let $\phi$ be the matrix satisfying
$\phi e_k=a_ke_k$ for $k=1,\ldots,n$, so that $E=\phi B^n$. Since
	\[
	\phi=(\det\phi)^{1/n}
	\frac{\phi}{(\det\phi)^{1/n}},
	\]
	the affine covariance and homogeneity of $S_{q-1}E$ \eqref{eqsec4-affine-property-of-S-q-minus-one-K} give
	\[
	S_{q-1}E
	=
	(\det\phi)^{1/(q-1)}\phi S_{q-1}B^n.
	\]
	By rotational symmetry, $S_{q-1}B^n=c_{q,n}B^n$ for some
	$c_{q,n}>0$. Since $\beta=(q-1)s$,
	\begin{align}
		\Psi_{s,q}(E)
		&=
		\frac{1}{n\omega_n}
		\int_{\Sn}\rho_{S_{q-1}E}(u)^\beta\,du \notag\\
		&=
		c_{q,n}^{\beta}(\det\phi)^s
		\frac{1}{n\omega_n}\int_{\Sn}\rho_E(u)^\beta\,du \notag\\
		&=
		c_{q,n}^{\beta}(a_1\cdots a_n)^s
		\widetilde V_\beta(E).
		\label{eqsec4-Psi-to-dual-quermassint}
	\end{align}
	
	Set
	\[
	D=
	\left\{
	x\in\operatorname{span}\{e_1,\ldots,e_{m-1}\}:
	\sum_{k=1}^{m-1}\frac{(x\cdot e_k)^2}{a_k^2}\le1
	\right\},
	\]
	\[
	J=[-a_m,a_m]e_m,
	\]
	and
	\[
	Y=
	\left\{
	y\in\operatorname{span}\{e_{m+1},\ldots,e_n\}:
	\sum_{k=m+1}^{n}\frac{(y\cdot e_k)^2}{a_k^2}\le1
	\right\}.
	\]
	If $m=1$, we set $D=\{o\}\subset\mathbb R^0$, and if $m=n$, we set
	$Y=\{o\}\subset\mathbb R^0$. We use the conventions
	$\vol_0(\{o\})=1$
	and that empty products are equal to $1$.
	
	Clearly,
	\[
	E\subset D\times J\times Y.
	\]
	Since $\beta>0$, \eqref{eqsec5-dual-quermassintegral-by-integration-on-K} gives
	\begin{align}
		\widetilde V_\beta(E)
		&\le
		\frac{\beta}{n\omega_n}
		\int_D\int_J\int_Y
		|(x,t,y)|^{\beta-n}\,dy\,dt\,dx \notag\\
		&\le
		\frac{\beta}{n\omega_n}
		\vol_{m-1}(D)
		\int_J\int_Y|(t,y)|^{\beta-n}\,dy\,dt,
		\label{eqsec4-dual-msint-without-D}
	\end{align}
	where the last line follows from $\beta-n<0$ and $|(x,t,y)|\ge |(t,y)|$. Moreover,
	\[
	\vol_{m-1}(D)
	=
	\omega_{m-1}\prod_{k=1}^{m-1}a_k.
	\]
	
	Suppose first that $m\le n-1$. Since $a_k\le1$, we have
	$Y\subset B^{n-m}$, and
	\[
	|(t,y)|\ge\frac{|t|+|y|}{2}.
	\]
	Using polar coordinates in $\operatorname{span}\{e_{m+1},\ldots,e_n\}$,
	we obtain
	\begin{align*}
		\int_J\int_Y|(t,y)|^{\beta-n}\,dy\,dt
		&\le
		2^{n-\beta}(n-m)\omega_{n-m}
		\int_{-a_m}^{a_m}\int_0^1
		(r+|t|)^{\beta-n}r^{n-m-1}\,dr\,dt.
	\end{align*}
	Since $r\le r+|t|$,
	\[
	(r+|t|)^{\beta-n}r^{n-m-1}
	\le
	(r+|t|)^{\beta-m-1}.
	\]
	As $m-1<\beta<m$,
	\begin{align*}
		\int_0^1(r+|t|)^{\beta-m-1}\,dr
		&=
		\frac{|t|^{\beta-m}-(1+|t|)^{\beta-m}}{m-\beta}\\
		&\le
		\frac{|t|^{\beta-m}}{m-\beta}.
	\end{align*}
	Therefore
	\begin{align}
		\int_J\int_Y|(t,y)|^{\beta-n}\,dy\,dt
		&\le
		\frac{(n-m)\omega_{n-m}2^{n-\beta}}
		{m-\beta}
		\int_{-a_m}^{a_m}|t|^{\beta-m}\,dt \notag\\
		&=
		\frac{(n-m)\omega_{n-m}2^{n-\beta+1}}
		{(m-\beta)(\beta+1-m)}
		a_m^{\beta+1-m}.
		\label{eqsec4-upperbound-for-JY}
	\end{align}
	Combining this with \eqref{eqsec4-dual-msint-without-D}, we get
	\[
	\widetilde V_\beta(E)
	\le
	C_{1}
	\left(\prod_{k=1}^{m-1}a_k\right)
	a_m^{\beta+1-m},
	\]where $$C_{1}=\frac{\beta\omega_{m-1}(n-m)\omega_{n-m}2^{n-\beta+1}}
	{n\omega_n(m-\beta)(\beta+1-m)}.$$
	
	If $m=n$, then $Y=\{o\}$ and
	\[
	\int_J |t|^{\beta-n}\,dt
	=
	\frac{2a_n^{\beta-n+1}}{\beta-n+1},
	\]
	so we have
	\[
	\widetilde V_\beta(E)
	\le
	C_{2}
	\left(\prod_{k=1}^{n-1}a_k\right)
	a_n^{\beta+1-n},
	\]where $$C_2=\frac{2\beta\omega_{n-1}}{n\omega_n(\beta-n+1)}.$$ 
Thus, in all cases, there exists a constant $C_{\beta,n}>0$,
depending only on $\beta$ and $n$, such that
\[
\widetilde V_\beta(E)
\le
C_{\beta,n}
\left(\prod_{k=1}^{m-1}a_k\right)
a_m^{\beta+1-m}.
\]
	Substituting this into
	\eqref{eqsec4-Psi-to-dual-quermassint}, we obtain
	\begin{align*}
		\Psi_{s,q}(E)
		&\le
		c_{s,q,n}
		(a_1\cdots a_n)^s
		\left(\prod_{k=1}^{m-1}a_k\right)
		a_m^{\beta+1-m}\\
		&=
		c_{s,q,n}
		\left(\prod_{k=1}^{m-1}a_k\right)^{s+1}
		a_m^{\beta+s+1-m}
		\left(\prod_{k=m+1}^{n}a_k\right)^s,
	\end{align*}where $c_{s,q,n}=c_{q,n}^\beta C_{\beta,n}$ depends only on $s,q,n$, since $\beta=(q-1)s$.
	This proves the lemma.
\end{proof}

For $0<s_1\le s_2$, H\"older's inequality gives
\begin{equation}\label{eqsec4-monotonicity-of-Psi-with-s}
\Psi_{s_1,q}(K)
\le
\Psi_{s_2,q}(K)^{s_1/s_2}.
\end{equation}

For a nonzero finite Borel measure $\mu$, define $L_{s,q,\mu}:\Kne\to\bR$ by
\begin{equation*}
L_{s,q,\mu}(K)
=
E_\mu(K)+\frac{1}{(n+q-1)s}\log\Psi_{s,q}(K).
\end{equation*}
This functional is $0$-homogeneous and is continuous in the Hausdorff metric.
For $K\in\Kno$, let
\[
R_K=\max\{|x|:x\in K\},
\qquad
r_K=\max\{r>0:rB^n\subset K\}.
\]

\begin{lemma}\label{lemsec4-degeneration}
Let $q\ge1$ and $s>0$, and assume that $(q-1)s\in[0,n]$. Let $\mu$ be a nonzero even finite Borel measure satisfying \eqref{eqsec4-subspace-mass-inequality}. If $(K_l)\subset\Kne$ is a sequence satisfying $R_{K_l}=1$ and $r_{K_l}\to0$, then
\[
\liminf_{l\to\infty} L_{s,q,\mu}(K_l)=-\infty.
\]
\end{lemma}

\begin{proof}
 Let $E_l$ be the John ellipsoid of $K_l$. Since $K_l$ is origin-symmetric, John's theorem gives
\begin{equation}\label{eqsec4-john-theorem}
E_l\subset K_l\subset\sqrt n\,E_l;
\end{equation}
see \cite{schneiderbook}. Write
\[
E_l
=
\left\{x\in\Rn:
\sum_{k=1}^n\frac{(x\cdot e_{kl})^2}{a_{k,l}^2}\le1
\right\},
\]
where $0<a_{1,l}\le\cdots\le a_{n,l}$ and
$(e_{1l},\ldots,e_{nl})$ is an ordered orthonormal basis. After passing to a subsequence and relabelling both
$(K_l)$ and $(E_l)$, we may assume that
$(e_{1l},\ldots,e_{nl})$ converges to an orthonormal
basis $(e_1,\ldots,e_n)$. From \eqref{eqsec4-john-theorem},
\[
\frac1{\sqrt n}\le a_{n,l}\le1,
\qquad
a_{1,l}\longrightarrow0.
\]
Since $E_\mu$ decreases and $\Psi_{s,q}$ increases under inclusion,
\begin{equation}\label{eqsec4-upperbound-for-GKl}
L_{s,q,\mu}(K_l)
\le
E_\mu(E_l)
+
\frac{1}{(n+q-1)s}\log\Psi_{s,q}(E_l)
+
\frac12\log n.
\end{equation}

 Let $\beta=(q-1)s$. First suppose that $\beta\notin\mathbb Z$ and put
$m=\lceil\beta\rceil$.
By Lemma~\ref{lemsec4-entropy-estimate}, there exist
$t_1>0$ and $c_1\in\mathbb R$, independent of $l$, such that,
for all sufficiently large $l$,
\begin{align*}
	E_\mu(E_l)
	\le{}&
	-\frac{s+1}{(n+q-1)s}
	\log\left(\prod_{k=1}^{m-1}a_{k,l}\right)-\frac{\beta+s+1-m}{(n+q-1)s}\log a_{m,l}\\
	&-\frac{s}{(n+q-1)s}
	\log\left(\prod_{k=m+1}^{n}a_{k,l}\right)
	+t_1\log a_{1,l}+c_1.
\end{align*}
On the other hand, Lemma~\ref{lemsec4-ellipsoid-estimate}
gives
\begin{align*}
	\frac{1}{(n+q-1)s}\log\Psi_{s,q}(E_l)
	\le{}&
	\frac{s+1}{(n+q-1)s}
	\log\left(\prod_{k=1}^{m-1}a_{k,l}\right)+\frac{\beta+s+1-m}{(n+q-1)s}\log a_{m,l}\\
	&+\frac{s}{(n+q-1)s}
	\log\left(\prod_{k=m+1}^{n}a_{k,l}\right)
	+\frac{\log c_{s,q,n}}{(n+q-1)s}.
\end{align*}
Adding these two inequalities, we obtain
\[
E_\mu(E_l)
+\frac{1}{(n+q-1)s}\log\Psi_{s,q}(E_l)
\le t_1\log a_{1,l}+c_1+\frac{\log c_{s,q,n}}{(n+q-1)s}.
\]
 Hence by the fact that $a_{1l}\to0$ and \eqref{eqsec4-upperbound-for-GKl}, \[
\liminf_{l\to\infty} L_{s,q,\mu}(K_l)=-\infty.
\]

Now suppose that $\beta\in\{1,\ldots,n-1\}$.
By Lemma~\ref{lemsec4-entropy-estimate}, with $m=\beta$,
there exist $t_2>0$ and $c_2\in\mathbb R$, independent of $l$,
such that, for all sufficiently large $l$,
\[
E_\mu(E_l)
\le
-\frac{s+1}{(n+q-1)s}\log\left(\prod_{k=1}^{m}a_{k,l}\right)
-\frac1{n+q-1}\log\left(\prod_{k=m+1}^{n}a_{k,l}\right)
+t_2\log a_{1,l}+c_2.
\] Now we set $\beta_1=(q-1)s_1$, where $s_1>s$ is chosen so that $m<\beta_1<m+1$. 
By \eqref{eqsec4-monotonicity-of-Psi-with-s} and
Lemma~\ref{lemsec4-ellipsoid-estimate} applied to $\Psi_{s_1,q}(E_l)$,
we obtain
\begin{align*}
\frac{1}{(n+q-1)s}\log\Psi_{s,q}(E_l)
\le{}&
\frac{s_1+1}{(n+q-1)s_1}
\log\left(\prod_{k=1}^{m}a_{k,l}\right)+
\frac{\beta_1+s_1-m}{(n+q-1)s_1}\log a_{m+1,l}\\
&
+
\frac{1}{n+q-1}\log\left(\prod_{k=m+2}^{n}a_{k,l}\right)
+C_1,
\end{align*}where $C_1$ is a constant independent of $l$.
Adding the two preceding estimates, we obtain
\begin{align*}
	E_\mu(E_l)
	+\frac{1}{(n+q-1)s}\log\Psi_{s,q}(E_l)
	\le{}&
	\frac{1}{n+q-1}
	\left(\frac1{s_1}-\frac1s\right)
	\sum_{k=1}^{m}\log a_{k,l}\\
	&+
	\frac{\beta_1-m}{(n+q-1)s_1}\log a_{m+1,l}
	+t_2\log a_{1,l}+C_3,
\end{align*}where $C_3$ is a constant independent of $l$.
Since $\beta_1>m$ and $a_{m+1,l}\le1$,
\[
\frac{\beta_1-m}{(n+q-1)s_1}\log a_{m+1,l}\le0,
\]
so this term may be discarded. Moreover, $s_1>s$ implies
\[
\frac1{s_1}-\frac1s<0.
\]
Since
\[
\log a_{1,l}\le\log a_{k,l},
\qquad k=1,\ldots,m,
\]
we have
\[
\left(\frac1{s_1}-\frac1s\right)
\sum_{k=1}^{m}\log a_{k,l}
\le
m\left(\frac1{s_1}-\frac1s\right)\log a_{1,l}.
\]
Consequently,
\begin{equation}\label{eqsec4 integer Lsqmu estimate}
	E_\mu(E_l)
	+
	\frac{1}{(n+q-1)s}\log\Psi_{s,q}(E_l)
	\le
	\left[
	t_2+\left(\frac1{s_1}-\frac1s\right)\frac{m}{n+q-1}
	\right]\log a_{1,l}+C_3.
\end{equation}We now choose $s_1$ sufficiently close to $s$ so that
\[
t_2+\left(\frac1{s_1}-\frac1s\right)\frac{m}{n+q-1}>0.
\]
 So the right-hand side of \eqref{eqsec4 integer Lsqmu estimate} tends to $-\infty$. \eqref{eqsec4-upperbound-for-GKl} therefore gives \[
\liminf_{l\to\infty} L_{s,q,\mu}(K_l)=-\infty.
\]

Finally, suppose that $(q-1)s\in\{0,n\}$.
In both cases, Lemma~\ref{lemsec4-entropy-estimate}
gives constants $t_0>0$ and $c_0\in\mathbb R$,
independent of $l$, such that, for all sufficiently large $l$,
\begin{equation}\label{eqsec4 0 n ellispoid entropy estimate}
	E_\mu(E_l)
	\le
	-\frac1n\log\left(\prod_{k=1}^n a_{k,l}\right)
	+t_0\log a_{1,l}+c_0.
\end{equation}
If $(q-1)s=0$, then $q=1$ and
\[
\Psi_{s,1}(E_l)
=V(E_l)^s
=\omega_n^s\left(\prod_{k=1}^n a_{k,l}\right)^s.
\]
If $(q-1)s=n$, by \eqref{eqsec4-Psi-to-dual-quermassint}, we have
\begin{align*}
\Psi_{s,q}(E_l)
&=c_{q,n}^{\,n}\left(\prod_{k=1}^n a_{k,l}\right)^s
\frac1{n\omega_n}\int_{\Sn}\rho_{E_l}(u)^n\,du\\
&=c_{q,n}^{\,n}\left(\prod_{k=1}^n a_{k,l}\right)^{s+1}.
\end{align*}
 Since $(n+q-1)s=ns$ when $(q-1)s=0$ and $(n+q-1)s=n(s+1)$ when $(q-1)s=n$, in either case
\begin{equation}\label{eqsec4 0 n elliposid estimate}
	\frac1{(n+q-1)s}\log\Psi_{s,q}(E_l)
	=\frac1n\log\left(\prod_{k=1}^n a_{k,l}\right)+C,
\end{equation}
where $C$ is a constant independent of $l$. Combining \eqref{eqsec4 0 n ellispoid entropy estimate}
and \eqref{eqsec4 0 n elliposid estimate}
with \eqref{eqsec4-upperbound-for-GKl}, we obtain
\[
L_{s,q,\mu}(K_l)
\le
t_0\log a_{1,l}+C+c_0+\frac12\log n
\longrightarrow-\infty
\qquad\text{as }l\to\infty.
\]
This completes the proof.
\end{proof}

\begin{theorem}\label{thmsec4-existence-of-maximizer}
Let $q\ge1$ and $s>0$ satisfy $(q-1)s\in[0,n]$. If $\mu$ is a nonzero even finite Borel measure satisfying \eqref{eqsec4-subspace-mass-inequality}, then there exists $K_0\in\Kne$ such that
\[
L_{s,q,\mu}(K_0)
=
\sup\{L_{s,q,\mu}(K):K\in\Kne\}.
\]
\end{theorem}

\begin{proof}
Let $(K_l)$ be a maximizing sequence of $L_{s,q,\mu}$. Since $E_\mu(B^n)=0$, the sequence $(L_{s,q,\mu}(K_l))$ is bounded below. By the $0$-homogeneity of $L_{s,q,\mu}$, we can normalize $(K_l)$ so that $R_{K_l}=1$. We claim that
$r_{K_l}$ is bounded below by a positive constant.
Otherwise, there would be a subsequence, again denoted
by $(K_l)$, such that $r_{K_l}\to0$.
Lemma~\ref{lemsec4-degeneration} would then imply
\[
\liminf_{l\to\infty}L_{s,q,\mu}(K_l)=-\infty,
\]
contradicting the fact that the sequence
$(L_{s,q,\mu}(K_l))$ is bounded from below. By the Blaschke selection theorem, a further subsequence converges to a body $K_0\in\Kne$. The Hausdorff continuity of $L_{s,q,\mu}$ shows that\[
L_{s,q,\mu}(K_0)
=
\sup\{L_{s,q,\mu}(K):K\in\Kne\}.
\]This completes the proof.
\end{proof}

\begin{proof}[Proof of Theorem~\ref{thmsec4-solution-to-log-mp}]
By Theorem~\ref{thmsec4-existence-of-maximizer}, let $K_0\in\Kne$ be a maximizer of $L_{s,q,\mu}$. For an even function $g\in C(\Sn)$, define
\[
h_t=h_{K_0}e^{tg},
\qquad
K_t=[h_t].
\]
For sufficiently small $|t|$, $h_t$ is positive and even, and $K_t\in\Kne$. Set
\[
\varphi(t)
=-\frac1{|\mu|}\int_{\Sn}\log h_t\,d\mu
+\frac{1}{(n+q-1)s}\log\Psi_{s,q}(K_t).
\]
Since $h_{K_t}\le h_t$,
\[
\varphi(t)
\le
L_{s,q,\mu}(K_t)
\le
L_{s,q,\mu}(K_0)
=
\varphi(0).
\]
Thus $\varphi'(0)=0$. By Theorem~\ref{thmsec2-variational-formula-for-Lp}, with $p=0$, we have
\[
-\frac1{|\mu|}\int_{\Sn}g(v)\,d\mu(v)
+
\frac{1}{(n+q-1)\Psi_{s,q}(K_0)}
\int_{\Sn}g(v)h_{K_0}(v)\,dF_{s,q}(K_0,v)
=0.
\]
This identity holds for every even continuous $g$. Both $\mu$ and the measure $h_{K_0}F_{s,q}(K_0,\cdot)$ are even, the latter by Proposition~\ref{prop-evenness-of-sq-chord-measure}. Hence
\[
\frac{|\mu|}{(n+q-1)\Psi_{s,q}(K_0)}
 h_{K_0}(v)\,dF_{s,q}(K_0,v)
=
d\mu(v).
\]
Choose $\lambda>0$ so that
\[
\lambda^{(n+q-1)s}
=
\frac{|\mu|}{(n+q-1)\Psi_{s,q}(K_0)}.
\]
By Proposition~\ref{prop-homogeneity-of-sq-chord-measure}, $K=\lambda K_0$ satisfies
\[
d\mu(v)=h_K(v)\,dF_{s,q}(K,v).
\] This proves the theorem.
\end{proof}

\section{\texorpdfstring{The asymptotic behavior of $(s,q)$-chord measures}{The asymptotic behavior of (s,q)-chord measures}}\label{sec:asymptotic}
In this section, we study the asymptotic behavior of the $(s,q)$-chord measure as $q\to0^+$. We first recall mixed volumes and mixed area measures. If $K_1,\ldots,K_m$ are compact convex sets and $\lambda_1,\ldots,\lambda_m\ge0$, then
\[
V(\lambda_1K_1+\cdots+\lambda_mK_m)
=
\sum_{i_1,\ldots,i_n=1}^{m}
\lambda_{i_1}\cdots\lambda_{i_n}
V(K_{i_1},\ldots,K_{i_n}),
\]
 where $V(K_{i_1},\ldots,K_{i_n})$ is the mixed volume, which is nonnegative, symmetric, and
 multi-linear in each argument. The mixed volume has the following representation
	\begin{equation}\label{eqsec5-mixedvolume-integral-representation}
	V(K_1,\ldots,K_n)
	=
	\frac1n\int_{\Sn}h_{K_1}(v)\,dS(K_2,\ldots,K_n;v),
\end{equation}where $S(K_2,\ldots,K_n;\cdot)$ is the mixed area measure.
We write
\[
S_{n-2}(K,L;\cdot)
=S(\underbrace{K,\ldots,K}_{n-2\ \mathrm{times}},L;\cdot),
\]and the $n-2$-th area measure $S_{n-2}(K,\cdot)$ is defined by\[
S_{n-2}(K,\cdot)=S_{n-2}(K,B^n;\cdot).\]
By multilinearity, mixed volumes and mixed area measures extend to the linear span of support functions. Thus, when $f_1,\ldots,f_n$ are differences of support functions, we can define $V(f_1,\ldots,f_n)$ and $S(f_1,\ldots,f_{n-1};\cdot)$; see \cite[Section~5.2]{schneiderbook}.
For later use, we also recall mixed discriminants. Let $A_1,\ldots,A_m$ be symmetric $d\times d$ matrices and $\lambda_1,\ldots,\lambda_m\ge0$. The mixed discriminant $\operatorname D$ is defined by
\begin{equation}\label{defsec5-mixed-discriminant}
	\det(\lambda_1A_1+\cdots+\lambda_mA_m)
	=
	\sum_{i_1,\ldots,i_d=1}^{m}
	\lambda_{i_1}\cdots\lambda_{i_d}
	\operatorname D(A_{i_1},\ldots,A_{i_d}).
\end{equation}
For $f\in C^2(\Sn)$, let\[
\overline f(x)=|x|\,f\!\left(\frac{x}{|x|}\right),
\qquad x\in\mathbb R^n\setminus\{0\}.
\] For $v\in\Sn$, $D^2\overline f(v)|_{v^\perp}$ denotes the restriction of the Euclidean Hessian to $v^\perp$. If $K_2,\ldots,K_n$ have $C^2$ boundaries and positive Gauss curvature, then
\begin{equation}\label{eqsec5-mixed-surface-area-measure-representation-by-mixed-discriminant}
	dS(K_2,\ldots,K_n;v)
	=
	\operatorname D\bigl(D^2\overline{h_{K_2}}(v)|_{v^\perp},\ldots,D^2\overline{h_{K_n}}(v)|_{v^\perp}\bigr)\,dv;
\end{equation}
see \cite{schneiderbook}.  

Throughout this section, fix $s\in\bR\setminus\{0\}$. We introduce the origin-symmetric zonoid $\bP_sK$ which is defined by
\begin{equation}\label{defsec5-s-zonoid}
	h_{\bP_sK}(v)
	=
	\int_{\Sn}|v\cdot u|h_{\Pi K}(u)^{s-1}\,du,
	\qquad v\in\Sn.
\end{equation}

\begin{definition}\label{defsec5-s-area-measure}
Let $K\in\Kn$. The $s$-area measure of $K$ is defined by
\[
S^s_{n-2}(K,\omega)
=
\frac{n-1}{2n\omega_n}
S_{n-2}(K,\bP_sK;\omega)
\]
for every Borel set $\omega\subset\Sn$.
\end{definition}

When $s=1$, $\bP_1K=2\omega_{n-1}B^n$, and hence
\[
S^1_{n-2}(K,\cdot)
=
\frac{(n-1)\omega_{n-1}}{n\omega_n}S_{n-2}(K,\cdot).
\]
Xi--Zhao \cite{XiZhao26} introduced the affine area measure,
which, in our notation, is defined by
\[
S^{(a)}_{n-2}(K,\cdot)
=
\frac{n-1}{4n}
S_{n-2}(K,\bP_{-n}K;\cdot).
\]
Therefore, when $s=-n$,
\[
S^{-n}_{n-2}(K,\cdot)
=
\frac{2}{\omega_n}S^{(a)}_{n-2}(K,\cdot).
\]
Thus the $s$-area measure also recovers the affine area
measure up to a multiplicative constant. Our main result is the following.
\begin{theorem}\label{thmsec5-limit-of-Fsq}
Let $K\in\Kn$ have $C^2$ boundary and everywhere positive Gauss curvature. Then, for every Borel set $\omega\subset\Sn$,
\[
\lim_{q\to0^+}F_{s,q}(K,\omega)
=
S^s_{n-2}(K,\omega).
\]
\end{theorem}

Under suitable regularity assumptions, for $z\in\partial K$,
let $\kappa_K(z)$ denote the Gauss curvature of $\partial K$
at $z$. If $u$ is a unit tangent vector to $\partial K$ at $z$,
let $\kappa_K(z,u)$ denote the normal curvature
in the direction $u$. 

\begin{corollary}[Corollary~4.2, \cite{XiZhao26}]\label{lemsec5-mixed-area-measure-of-interval}
Suppose $K\in\Kn$ has $C^2$ boundary and everywhere positive Gauss curvature. Let $u\in\Sn$. For every Borel set $\omega\subset\Sn$,
\[
S_{n-2}(K,[-u,u];\omega)
=
\frac{2}{n-1}
\int_{\omega\cap u^\perp}
\frac{\kappa_K(z(v),u)}{\kappa_K(z(v))}\,d\cH^{n-2}(v),
\]
where $z(v)\in\partial K$ is the unique point with outer unit normal $v$.
\end{corollary}
Let
\[
\mathcal I=\{(u,v)\in\Sn\times\Sn:u\cdot v=0\}.
\]For every integrable function $H$ on $\mathcal I$, we have

\begin{equation}\label{eq-incidence-Fubini}
	\int_{\Sn}\int_{\Sn\cap u^\perp}H(u,v)\,d\cH^{n-2}(v)\,du
	=
	\int_{\Sn}\int_{\Sn\cap v^\perp}H(u,v)\,d\cH^{n-2}(u)\,dv,
\end{equation}
see~\cite[Theorem~7.1.1]{SchneiderWeil2008}. 
\begin{theorem}\label{thmsec5-s-area-measure-representation-of-curvature}
Suppose $K\in\Kn$ has $C^2$ boundary and everywhere positive Gauss curvature. Then, for every Borel set $\omega\subset\Sn$,
\[
S^s_{n-2}(K,\omega)
=
\frac{1}{n\omega_n}
\int_{\nu_K^{-1}(\omega)}
\int_{\Sn\cap\nu_K(z)^\perp}
h_{\Pi K}(u)^{s-1}\kappa_K(z,u)\,d\cH^{n-2}(u)\,d\cHnmone(z).
\]
\end{theorem}

\begin{proof}
Writing $K[n-2]$ for $n-2$ copies of $K$, Definition~\ref{defsec5-s-area-measure}, \eqref{eqsec5-mixedvolume-integral-representation}, 
\eqref{defsec5-s-zonoid}, and the symmetry of mixed volumes give, for every convex body $L$,
\begin{align*}
\int_{\Sn} h_L(v)\,dS^s_{n-2}(K,v)
&= \frac{n-1}{2\omega_n}
V(L,K[n-2],\bP_sK) \\
&= \frac{n-1}{2\omega_n}
V(\bP_sK,K[n-2],L) \\
&= \frac{n-1}{2n\omega_n}
\int_{\Sn} h_{\bP_sK}(v)\,dS_{n-2}(K,L;v) \\
&= \frac{n-1}{2n\omega_n}
\int_{\Sn} h_{\Pi K}(u)^{s-1}
\int_{\Sn} |v\cdot u|\,dS_{n-2}(K,L;v)\,du \\
&= \frac{n-1}{2\omega_n}
\int_{\Sn} h_{\Pi K}(u)^{s-1}
V([-u,u],K[n-2],L)\,du \\
&= \frac{n-1}{2\omega_n}
\int_{\Sn} h_{\Pi K}(u)^{s-1}
V(L,K[n-2],[-u,u])\,du \\
&= \frac{n-1}{2n\omega_n}
\int_{\Sn} h_{\Pi K}(u)^{s-1}
\int_{\Sn} h_L(v)\,dS_{n-2}(K,[-u,u];v)\,du.
\end{align*}
Since the linear span of support functions is dense in $C(\Sn)$
 (see \cite{schneiderbook}), the preceding identity extends to every
$\varphi\in C(\Sn)$,
\[
\int_{\Sn}\varphi(v)\,dS^s_{n-2}(K,v)
=
\frac{n-1}{2n\omega_n}
\int_{\Sn}h_{\Pi K}(u)^{s-1}
\int_{\Sn}\varphi(v)\,dS_{n-2}(K,[-u,u];v)\,du.
\]
Consequently, for every Borel set $\omega\subset\Sn$,
\[
S^s_{n-2}(K,\omega)
=
\frac{n-1}{2n\omega_n}
\int_{\Sn}S_{n-2}(K,[-u,u];\omega)\,
h_{\Pi K}(u)^{s-1}\,du.
\]
Using Corollary~\ref{lemsec5-mixed-area-measure-of-interval} and \eqref{eq-incidence-Fubini},
\begin{align*}
S^s_{n-2}(K,\omega)
&=
\frac{1}{n\omega_n}
\int_{\Sn}\int_{\Sn\cap u^\perp}
\mathbf 1_\omega(v)
\frac{\kappa_K(z(v),u)}{\kappa_K(z(v))}
\,d\cH^{n-2}(v)
h_{\Pi K}(u)^{s-1}\,du\\
&=
\frac{1}{n\omega_n}
\int_{\Sn}\mathbf 1_\omega(v)
\frac{1}{\kappa_K(z(v))}
\int_{\Sn\cap v^\perp}
\kappa_K(z(v),u)h_{\Pi K}(u)^{s-1}\,d\cH^{n-2}(u)\,dv.
\end{align*}
Finally, using the change of  variables $v=\nu_K(z)$, we complete the proof.
\end{proof}

Under suitable regularity assumptions, we have the following variational formula for $\Psi_{s,0}$.

\begin{theorem}\label{thmsec5-variation-at-zero}
Suppose $K\in\Kn$ has $C^2$ boundary and everywhere positive Gauss curvature. Let $g\in C^2(\Sn)$, $h_t=h_K+tg$, and $K_t=[h_t]$. Then,
\[
\left.\frac{d}{dt}\right|_{t=0}\Psi_{s,0}(K_t)
=
s\int_{\Sn}g(v)\,dS^s_{n-2}(K,v).
\]
\end{theorem}

\begin{proof}
Since $K$ has $C^2$ boundary and positive Gauss curvature, $h_t$ is the support function of $K_t$ for sufficiently small $|t|$. The Cauchy projection formula gives
\begin{equation}\label{eqsec6 n-1 dim V of Kt|u perp}
	\vol_{n-1}(K_t|u^\perp)
	=
	\frac12\int_{\Sn}|u\cdot v|
	\det\bigl(D^2\overline{h_t}(v)|_{v^\perp}\bigr)\,dv;
\end{equation}
see \cite{schneiderbook}. For a matrix or a convex body $A$ and an integer $r\ge0$, we use $A[r]$ to denote $r$ copies of $A$. Differentiating \eqref{eqsec6 n-1 dim V of Kt|u perp} and using \eqref{defsec5-mixed-discriminant} and \eqref{eqsec5-mixedvolume-integral-representation},
\begin{align*}
\left.\frac{d}{dt}\right|_{t=0}\vol_{n-1}(K_t|u^\perp)
&=
\frac{n-1}{2}
\int_{\Sn}|u\cdot v|
\operatorname D\bigl(
D^2\overline{h_K}(v)|_{v^\perp}[n-2],
D^2\overline g(v)|_{v^\perp}
\bigr)\,dv\\
&=
\frac{n(n-1)}2V([-u,u],K[n-2],g)\\
&=
\frac{n-1}{2}
\int_{\Sn}g(v)\,dS_{n-2}(K,[-u,u];v)\\
&=
\int_{\Sn\cap u^\perp}
g(v)\frac{\kappa_K(z(v),u)}{\kappa_K(z(v))}\,d\cH^{n-2}(v).
\end{align*}

Differentiating $\Psi_{s,0}$ and using \eqref{eq-incidence-Fubini} now yields
\begin{align*}
\left.\frac{d}{dt}\right|_{t=0}\Psi_{s,0}(K_t)
&=
\frac{s}{n\omega_n}
\int_{\Sn}\vol_{n-1}(K|u^\perp)^{s-1}
\int_{\Sn\cap u^\perp}
g(v)
\frac{\kappa_K(z(v),u)}{\kappa_K(z(v))}\,d\cH^{n-2}(v)\,du\\
&=
\frac{s}{n\omega_n}
\int_{\Sn}\frac{g(v)}{\kappa_K(z(v))}
\int_{\Sn\cap v^\perp}
\vol_{n-1}(K|u^\perp)^{s-1}
\kappa_K(z(v),u)\,d\cH^{n-2}(u)\,dv\\
&=
\frac{s}{n\omega_n}
\int_{\partial K}g(\nu_K(z))
\int_{\Sn\cap\nu_K(z)^\perp}
h_{\Pi K}(u)^{s-1}\kappa_K(z,u)\,d\cH^{n-2}(u)\,d\cHnmone(z)\\
&=
s\int_{\Sn}g(v)\,dS^s_{n-2}(K,v),
\end{align*}
where the last equality follows from Theorem~\ref{thmsec5-s-area-measure-representation-of-curvature}.
\end{proof}

The following asymptotic lemma is the key step in the proof of Theorem~\ref{thmsec5-limit-of-Fsq}.

\begin{lemma}\label{lemsec5-lim-of-mixed-S-q-minus-one-and-K}
Let $K\in\Kn$ have $C^2$ boundary and everywhere positive Gauss curvature. Then, for every $z\in\partial K$,
\[
\lim_{q\to0^+}
q\int_{\Sn}
I_q(K,u)^{s-1}X_K(z,u)^{q-1}\,du
=
\int_{\Sn\cap\nu_K(z)^\perp}
h_{\Pi K}(\theta)^{s-1}\kappa_K(z,\theta)\,d\cH^{n-2}(\theta).
\]
\end{lemma}

\begin{proof}
	Let $\epsilon>0$. Since the integrands involved are invariant under translation,
	replacing $K$ with $K-z$, we may assume $z=o$, where $o$ is the origin.
	As $K$ has $C^2$ boundary, there exists an orthonormal basis
	$(e_1,\ldots,e_n)$ such that $e_n=-\nu_K(o)$ and, in a neighborhood of $o$,
	the boundary $\partial K$ can be represented as the graph of the following function
	\begin{equation}\label{eqsec5-graph-representation-of-boundary}
		y_n=f(y')
		=
		\frac12(\kappa_1y_1^2+\cdots+\kappa_{n-1}y_{n-1}^2)
		+o(|y'|^2),
	\end{equation}
	where $y=(y',y_n)=(y_1,\ldots,y_{n-1},y_n)$ and $\kappa_i$ are the
	principal curvatures of $K$ at $o$.

	 Let $\lambda>0$ be such that $\operatorname{diam}(\lambda K)\le1$. Since $X_{\lambda K}(\cdot,\cdot)\le1$, the functions $I_q(\lambda K,\cdot)$ increase as \(q\downarrow0\) and the limit is $I_0(\lambda K,\cdot)=h_{\Pi \lambda K}$. Since \(I_q(\lambda K,\cdot)\) and \(I_0(\lambda K,\cdot)\) are continuous on the unit sphere, \(I_q(\lambda K,\cdot)\to I_0(\lambda K,\cdot)\) uniformly by Dini's theorem. By the homogeneity of the directional chord integral, we have \(I_q(K,\cdot)\to I_0(K,\cdot)=h_{\Pi K}\) uniformly as \(q\to0^+\). Consequently, for sufficiently small $q$, the functions $I_q(K,\cdot)$ are bounded from above and away from \(0\) uniformly. By this uniform convergence and the uniform bounds, we may take $q>0$ sufficiently small so that
	 \begin{equation}\label{eqsec5-uniform-convergence-of-Iq}
	 	\left|
	 	I_q(K,u)^{s-1}
	 	-
	 	h_{\Pi K}(u)^{s-1}
	 	\right|
	 	<\epsilon,
	 	\qquad u\in\Sn.
	 \end{equation}
	
	Since, for $u\in\Sn$,
	\[
	X_K(o,u)=\rho_{K}(u)
	\quad\text{if }u\cdot e_n>0,
	\]
	and
	\[
	X_K(o,u)=\rho_{K}(-u)
	\quad\text{if }u\cdot e_n<0,
	\]
	and $I_q(K,-u)=I_q(K,u)$, we have
	\begin{equation}\label{eqsec5-integral-of-X-and-rho}
		q\int_{\Sn}I_q(K,u)^{s-1}X_K(o,u)^{q-1}\,du
		=
		2q\int_{\{u\in\Sn:u\cdot e_n>0\}}
		I_q(K,u)^{s-1}\rho_{K}(u)^{q-1}\,du.
	\end{equation}
	
	In a neighborhood of $o$ where
	\eqref{eqsec5-graph-representation-of-boundary} represents $\partial K$,
	we have
	\[
	\nu_K(y)
	=
	\frac{(\nabla f(y'),-1)}
	{\sqrt{1+|\nabla f(y')|^2}},
	\qquad
	d\cHnmone(y)
	=
	\sqrt{1+|\nabla f(y')|^2}\,dy'.
	\]
	Let $y'=r\theta$, where $r>0$ and
	$\theta\in\mathbb S^{n-2}$.
	
	Since $f(r\theta)=O(r^2)$,
	\[
	\frac{(r\theta,f(r\theta))}
	{\sqrt{r^2+f(r\theta)^2}}=	\frac{(\theta,f(r\theta)/r)}
	{\sqrt{1+f(r\theta)^2/r^2}}
	\longrightarrow
	(\theta,0)
	\]
	uniformly in $\theta$ as $r\to0^+$. By the
	uniform continuity of $h_{\Pi K}^{s-1}$ on $\Sn$, there exists
	$\delta>0$ such that, whenever $r<\delta$,
	\begin{equation}\label{eqsec5-uniform-continuity-of-hPi}
		\left|
		h_{\Pi K}\left(
		\frac{(r\theta,f(r\theta))}
		{\sqrt{r^2+f(r\theta)^2}}
		\right)^{s-1}
		-
		h_{\Pi K}(\theta,0)^{s-1}
		\right|
		<\epsilon.
	\end{equation}
	
	In a neighborhood of $o$, we have
	\begin{equation}\label{eqsec5-relation-between-rho-and-f}
		\rho_{K}(u)u
		=
		y=(y',f(y')).
	\end{equation}
	Moreover,
	\[
	\rho_{K}(u)^n\,du
	=
	|y\cdot\nu_K(y)|\,d\cHnmone(y),
	\qquad
	y'\cdot\nabla f(y')
	=
	2f(y')+o(|y'|^2).
	\]
	Hence
	\begin{equation}\label{eqsec5-change-variable-between-rho-and-y}
		\rho_{K}(u)^n\,du
		=
		\bigl(f(y')+o(|y'|^2)\bigr)\,dy'.
	\end{equation}
	
	Let
	\begin{align*}
		D_1
		&=
		\left\{
		\frac{(r\theta,f(r\theta))}{|(r\theta,f(r\theta))|}:
		\theta\in\mathbb S^{n-2},\ 0<r<\delta
		\right\},\\
		D_2
		&=
		\{u\in\Sn:u\cdot e_n>0\}\setminus D_1.
	\end{align*}
	Here $\delta>0$ is sufficiently small that the chosen
	local graph is defined for $|y'|<\delta$ and contains
	$\partial K\cap B(o,\delta)$.
	If $0<\rho_K(u)<\delta$, then
	$\rho_K(u)u=(r\theta,f(r\theta))$ for some $0<r<\delta$,
	so $u\in D_1$.
	Consequently, $\rho_K(u)\ge\delta$ for every $u\in D_2$.
	
	We first estimate the integral on $D_1$. By
	\eqref{eqsec5-uniform-convergence-of-Iq} and
	\eqref{eqsec5-change-variable-between-rho-and-y},
	\begin{align*}
		&q\int_{D_1}
		I_q(K,u)^{s-1}\rho_{K}(u)^{q-1}\,du\\
		&\le
		q\int_{D_1}
		\bigl(h_{\Pi K}(u)^{s-1}+\epsilon\bigr)
		\rho_{K}(u)^{q-1}\,du\\
		&=
		q\int_{\mathbb S^{n-2}}\int_0^\delta
		\left(
		h_{\Pi K}\left(
		\frac{(r\theta,f(r\theta))}
		{|(r\theta,f(r\theta))|}
		\right)^{s-1}
		+\epsilon
		\right)\\
		&\qquad\qquad\times
		|(r\theta,f(r\theta))|^{q-1-n}
		\bigl(f(r\theta)+o(r^2)\bigr)
		r^{n-2}\,dr\,d\cH^{n-2}(\theta).
	\end{align*}
	
	Since
	\[
	\kappa_K\left(o,(\theta,0)\right)
	=
	\kappa_1((\theta,0)\cdot e_1)^2+\cdots+
	\kappa_{n-1}((\theta,0)\cdot e_{n-1})^2,
	\]
	we have
	\[
	f(r\theta)
	=
	\frac12r^2\kappa_K\left(o,(\theta,0)\right)+o(r^2).
	\]
	Therefore
	\begin{align*}
		&q\int_{D_1}
		I_q(K,u)^{s-1}\rho_{K}(u)^{q-1}\,du\\
		&\le
		q\int_{\mathbb S^{n-2}}\int_0^\delta
		\left(
		h_{\Pi K}\left(
		\frac{(r\theta,f(r\theta))}
		{|(r\theta,f(r\theta))|}
		\right)^{s-1}
		+\epsilon
		\right)\\
		&\qquad\times
		\left(
		1+\frac{r^2\kappa_K(o,(\theta,0))^2}{4}+o(r^2)
		\right)^{\frac{q-1-n}{2}}
		\left(
		\frac{\kappa_K(o,(\theta,0))}2+o(1)
		\right)
		r^{q-1}\,dr\,d\cH^{n-2}(\theta).
	\end{align*}
	Let $\epsilon\in(0,1)$ be sufficiently small. Since
	\[
	1+\frac{r^2\kappa_K(o,(\theta,0))^2}{4}+o(r^2)
	=
	1+\frac{f(r\theta)^2}{r^2}\ge1,
	\]
	we may choose $\delta>0$ such that
\begin{equation}\label{eqsec6 tow sided bound of (y',f(y'))}
	1-\epsilon
	\le
	\left(
	1+\frac{r^2\kappa_K(o,(\theta,0))^2}{4}+o(r^2)
	\right)^{\frac{q-1-n}{2}}
	\le1
\end{equation}
	for $0<r<\delta$ and $0<q\le1$, uniformly in $\theta$.
	Moreover,
\begin{equation}\label{eqsec6 tow sided bound of guass curvature}
	\frac{\kappa_K(o,(\theta,0))}2-\epsilon
	\le
	\frac{\kappa_K(o,(\theta,0))}2+o(1)
	\le
	\frac{\kappa_K(o,(\theta,0))}2+\epsilon
\end{equation}
	for $0<r<\delta$, uniformly in $\theta$.
	
	Using the upper bounds of \eqref{eqsec6 tow sided bound of (y',f(y'))} and \eqref{eqsec6 tow sided bound of guass curvature} together with
	\eqref{eqsec5-uniform-continuity-of-hPi}, we obtain, for
	$q>0$ sufficiently small,
	\begin{align*}
		&q\int_{D_1}
		I_q(K,u)^{s-1}\rho_K(u)^{q-1}\,du\\
		&\quad\le
		q\int_{\mathbb S^{n-2}}\int_0^\delta
		\bigl(h_{\Pi K}(\theta,0)^{s-1}+2\epsilon\bigr)
		\left(\frac{\kappa_K(o,(\theta,0))}2+\epsilon\right)
		r^{q-1}\,dr\,d\cH^{n-2}(\theta)\\
		&\quad=
		\delta^q
		\int_{\mathbb S^{n-2}}
		\bigl(h_{\Pi K}(\theta,0)^{s-1}+2\epsilon\bigr)
		\left(\frac{\kappa_K(o,(\theta,0))}2+\epsilon\right)
		\,d\cH^{n-2}(\theta).
	\end{align*}
	
	Since $\rho_K(u)\ge\delta>0$ for $u\in D_2$ and
	$I_q(K,\cdot)^{s-1}$ is uniformly bounded for sufficiently
	small $q$, we have
	\[
	\lim_{q\to0^+}
	q\int_{D_2}I_q(K,u)^{s-1}\rho_K(u)^{q-1}\,du=0.
	\]
 Letting
	$q\to0^+$, and then letting $\epsilon\to0^+$, we obtain
	\begin{equation}\label{eqsec5-limsup-bound}
		\limsup_{q\to0^+}
		q\int_{\{u\in\Sn:u\cdot e_n>0\}}
		I_q(K,u)^{s-1}\rho_K(u)^{q-1}\,du
		\le
		\frac12
		\int_{\mathbb S^{n-2}}
		h_{\Pi K}(\theta,0)^{s-1}\kappa_K(o,(\theta,0))
		\,d\cH^{n-2}(\theta).
	\end{equation}
	Similarly, using the corresponding lower bounds of \eqref{eqsec5-uniform-convergence-of-Iq},
	\eqref{eqsec5-uniform-continuity-of-hPi}, \eqref{eqsec6 tow sided bound of (y',f(y'))} and \eqref{eqsec6 tow sided bound of guass curvature},
	we obtain
\begin{equation}\label{eqsec5 liminf-bound}
	\liminf_{q\to0^+}
	q\int_{\{u\in\Sn:u\cdot e_n>0\}}
	I_q(K,u)^{s-1}\rho_K(u)^{q-1}\,du
	\ge
	\frac12
	\int_{\mathbb S^{n-2}}
	h_{\Pi K}(\theta,0)^{s-1}\kappa_K(o,(\theta,0))
	\,d\cH^{n-2}(\theta).
\end{equation}
	Combining \eqref{eqsec5-limsup-bound} and \eqref{eqsec5 liminf-bound} with
	\eqref{eqsec5-integral-of-X-and-rho}, we obtain
	\[
	\lim_{q\to0^+}
	q\int_{\Sn}
	I_q(K,u)^{s-1}X_K(o,u)^{q-1}\,du
	=
	\int_{\mathbb S^{n-2}}
	h_{\Pi K}(\theta,0)^{s-1}
	\kappa_K(o,(\theta,0))\,d\cH^{n-2}(\theta).
	\]
This completes the proof.
\end{proof}

We recall the generalized dominated convergence theorem: if $f_j\to f$ and $\phi_j\to\phi$ almost everywhere, $|f_j|\le\phi_j$, the functions $\phi_j$ and $\phi$ are integrable, and $\int\phi_j\to\int\phi$, then $f$ is integrable and $\int f_j\to\int f$.

\begin{proof}[Proof of Theorem~\ref{thmsec5-limit-of-Fsq}]
Since $K$ has nonempty interior, $I_q(K,\cdot)$ and $h_{\Pi K}(\cdot)$ are bounded above and below away from $0$. Together with the fact that $I_q(K,\cdot)\to h_{\Pi K}(\cdot)$ uniformly as $q\to0^+$, this implies that there exists $C>0$ such that for every $z\in\partial K$,
\begin{align*}
q\int_{\Sn}
I_q(K,u)^{s-1}X_K(z,u)^{q-1}\,du
&\le
Cq\int_{\Sn}X_K(z,u)^{q-1}\,du\\
&=
2n\omega_nCq\,\widetilde V_{q-1}(K,z).
\end{align*}
It was shown in \cite[Corollary~4.14]{XiCPAM} that \[\lim_{q\to0^+}\int_{\nu^{-1}_{K}(\omega)}q\widetilde{V}_{q-1}(K,z)\,d\cH^{n-1}(z)=\int_{\nu^{-1}_{K}(\omega)}\lim_{q\to0^+}q\widetilde{V}_{q-1}(K,z)\,d\cH^{n-1}(z).\] 
Therefore the generalized dominated convergence theorem, Lemma~\ref{lemsec5-lim-of-mixed-S-q-minus-one-and-K} and Theorem \ref{thmsec5-s-area-measure-representation-of-curvature} imply
\begin{align*}
\lim_{q\to0^+}F_{s,q}(K,\omega)
&=
\frac{1}{n\omega_n}
\int_{\nu_K^{-1}(\omega)}
\int_{\Sn\cap\nu_K(z)^\perp}
h_{\Pi K}(u)^{s-1}\kappa_K(z,u)\,d\cH^{n-2}(u)\,d\cHnmone(z)\\
&=
S^s_{n-2}(K,\omega).
\end{align*}
This completes the proof.
\end{proof}
\section*{AI disclosure}
During the preparation of this paper, AI tools were used only
for proof checking and grammar checking. All proofs were written
by the authors, who take full responsibility for the content
of the paper.

	\end{document}